\documentclass[a4,12pt,multicol,makeidx]{article}
\def\origin{%
\hbox{}\vskip-\baselineskip\vskip-\topskip%
  \vbox to 0pt{\vskip-1in%
    \hbox to 0pt{\hskip-1in%
      \hbox to 0pt{\vrule width 1cm height .4pt depth 0mm\hss}%
      \vbox to 0pt{\hrule width .4pt height 0pt depth 1cm\vss}%
    \hss}%
  \vss}
  \vskip-\baselineskip
  \vbox to 0pt{\vskip-1in\vskip3cm%
    \hbox to 0pt{\hskip-1in\hskip3cm%
      \hbox to 0pt{\hss\vrule width 2cm height .4pt depth 0mm\hss}%
      \vbox to 0pt{\vss\hrule width .4pt height 1cm depth 1cm\vss}%
    \hss}%
  \vss}%
  \vskip5mm\hskip10mm (3cm,3cm)
}%
\def\CA{\rm A} \def\a{\alpha} \def\d{\delta} \def\l{\lambda} \def\G{\Gamma} \def\L{\Lambda} \def\p{\partial} \def\g{\gamma} \def\e{\varepsilon} 
\def\n{\nabla} \def\H{\nabla^2} \def\Rn{{\bf R}^n} \def\R{{\bf R}} 
\def\le{\underline{<}}

\newenvironment{theorem}{%
\par \bigskip \it}{%
\bigskip \par}
\newenvironment{proposition}{%
\par \bigskip \it}{%
\bigskip \par}
\newenvironment{lemma}{%
\par \bigskip \it}{%
\bigskip \par}

\newenvironment{example}{%
\par \bigskip \it}{%
\bigskip \par}

\makeindex
\title{Long time averaged reflection force\\and homogenization of  oscillating\\
Neumann boundary conditions.}
\author{Mariko Arisawa\\ \bigskip\\ GSIS, Tohoku University
\\ Aramaki 09, Aoba-ku, Sendai 980-8579, JAPAN\\
E-mail: arisawa@math.is.tohoku.ac.jp}
\date{}
\begin{document}
\maketitle
\bigskip
{\bf Abstract.} This paper concerns with two issues. The first issue is the existence and 
the uniqueness of the ergodic type number $d$ which appears in the oblique  boundary 
condition. The second issue is the application of the number for the study of homogenizations of oscillating Neumann boundary conditions. \\

{\bf R\'{e}sum\'{e}.} Dans cette article, nous traitons deux probl\`{e}mes. Le premier est 
 l'existence et l'unicit\'{e} d'un nombre du type ergodique $d$ qui appara\^{i}t dans la condition oblique sur le bord. Le second est l'application de ce nombre pour la recherche 
 des homog\'{e}n\'{e}izationses conditions Neumann sur des bords oscillants. \\ 

\newpage
\section{Introduction}
	First, we are concerned with the existence and uniqueness of the 
number $d$ in the following problem.

\begin{equation}\label{1}
F(x,\n u, \H u)=0 \qquad \hbox{in} \quad \Omega,
\end{equation} 
\begin{equation}\label{2}
d+<\n u,\g(x)>-g(x)=0 \qquad \hbox{on} \quad \partial \Omega,
\end{equation}
where $\Omega$ is a domain in $\Rn$,  $F$ is a fully nonlinear uniformly elliptic Hamilton-Jacobi-Bellman (HJB in short) operator:
\begin{equation}\label{hjb}
F(x,\n u, \H u)=\sup_{\a\in \CA} \{
-\sum_{i,j=1}^n a_{ij}^{\a}(x)\frac{\p^2 u}{\p x_i\p x_j}-\sum_{i=1}^n b_i^{\a}(x)\frac{\p u}{\p x_i}\},
\end{equation}
satisfying the following conditions. 
 $\CA$ is a set of controls, and by denoting $n\times n$ matrices $A^{\a}=({a_{ij}^{\a}(x)})_{ij}$ ($\a\in \CA$), there exist $n\times m$ matrices 
$\sigma^{\a}$ such that 
\begin{eqnarray}\label{unif}
&&A^{\a}(x)=\sigma^{\a}(\sigma^{\a})^t (x) \qquad \hbox{any} \quad x\in \Omega,\quad \a\in 
\CA,\nonumber\\
&&\l_1 I\le A^{\a}(x)  \le \L_1 I \qquad \hbox{any} \quad x\in \Omega,\quad \a\in 
\CA,
\end{eqnarray}
where $0<\l_1\le \L_1$ positive constants, $I$ the $n\times n$ identity matrix. 
 There exists a positive constant $L>0$ such that 
\begin{eqnarray}\label{lips}
|a_{ij}^{\a}(x)-a_{ij}^{\a}(y)|&\le & L|x-y| \qquad \hbox{any}\quad 1\le i,j\le n,\quad x\in \Omega,\quad \a\in A,\nonumber\\
|b_{i}^{\a}(x)-b_{i}^{\a}(y)|&\le & L|x-y| \qquad \hbox{any}\quad 1\le i\le n,\quad x\in \Omega,\quad \a\in A.
\end{eqnarray}
 	There also exists a positive constant $\g_0$, such that for the outward unit normal vector ${\bf n}(x)$ ($x\in \p\Omega$), $\g(x)$ satisfies
\begin{equation}\label{gamma}
<\g(x),{\bf n}(x)> \quad \ge \g_0>0 \qquad \hbox{any}\quad x\in \p\Omega. 
\end{equation} 
The domain  $\Omega$ is assumed to be either one of the following: 
\begin{equation}\label{doma}
\hbox{Bounded open domain in} \quad \Rn \quad 
 \quad \hbox{with}\quad C^{3,1} \quad\hbox{boundary}, 
\end{equation}
or
\begin{eqnarray}\label{dom2}
&&\hbox{Half space in}\quad  \Rn, \quad \hbox{periodic in the first} \quad n-1 \quad \hbox{variables with}\quad C^{3,1}\nonumber \\
&&\hbox{boundary}\nonumber\\
&&: \{(x',x_n)|\quad\hbox{periodic in}\quad x'=(x_1,...,x_{n-1})\in (\R/{\bf Z})^{n-1}, \quad x_n\ge f_1(x')\}, \nonumber\\
&&\hbox{where}\quad f_1\in C^{3,1}((\R/{\bf Z})^{n-1})).
\end{eqnarray}

 (In the latter case (\ref{dom2}), a supplement boundary condition at $x_n=\infty$ will be added to (\ref{1})-(\ref{2}).) \\

	The following example implies the qualitative meaning of the number $d$.

\begin{example}{\bf Example 1.1.} Let $\Omega$ be a domain in $(\ref{doma})$, and $g(x)$ be a Lipschitz continuous function on $\p\Omega$. Assume that there exists a number $d$  such that  the following problem 
has a viscosity solution.
$$
	-\Delta u=0 \qquad \hbox{in} \quad \Omega,
$$
$$
	d+<\n u,{\bf n}(x)>-g(x)=0 \qquad \hbox{on}\quad \p\Omega.
$$
Then, 
$$
	d=\frac{1}{|\p\Omega|}\int_{\p\Omega} g(x) dS.
$$
\end{example}
{\it Proof of  Example  1.1.}  In the Green's first identity: 
$$
	\int_{\Omega} \Delta u v dx+\int_{\Omega} \n u\cdot\n v dx =
	\int_{\p\Omega}v\frac{\p u}{\p n} dS,
$$ 
we put $v=1$, and get $d|\p\Omega|=\int_{\p\Omega}g(x)dS$.\\

	Thus, $d$ is a kind of the averaged quantity on $\p\Omega$. For general Hamiltonians $F$, the way to construct 
 the number $d$ and $u(x)$ in (\ref{1})-(\ref{2}) is the following. Here we assume that (\ref{doma}) holds. (The case (\ref{dom2}) is more complicated, and will be treated  in   Section 3 below.)
 For any $\lambda>0$, consider 
\begin{equation}\label{le}
F(x,\n u_{\l}, \H u_{\l})=0 \qquad \hbox{in} \quad \Omega,
\end{equation} 
\begin{equation}\label{lb}
\l u_{\l}+<\n u,\g(x)>-g(x)=0 \qquad \hbox{on} \quad \partial \Omega.
\end{equation}
 The regularity of $u_{\l}$ ($\l\in(0,1)$) which will be shown in Section 2 yields, for any  fixed $x_0\in \Omega$ 
\begin{equation}\label{lim}
\lim_{\l\downarrow 0} \l u_{\l}(x)=d \qquad \hbox{uniformly in}\quad \overline{\Omega},
\end{equation}
and by taking a subsequence $\l'\downarrow 0$,
\begin{equation}\label{ulim}
\lim_{\l'\downarrow 0} (u_{\l'}(x)-u_{\l'}(x_0))=u(x) \qquad \hbox{uniformly in}\quad \overline{\Omega}.
\end{equation}
 The limit number $d$ is unique in the sense that with which  (\ref{1})-(\ref{2}) has a viscosity solution. The above limit function $u(x)$ is one of such  solutions. (The solution of (\ref{1})-(\ref{2}) is not unique, for $u+C$ ($C$ constant) is also a solution.)  We shall  show in Section 2 these facts.  Now, the meaning of the number $d$ can be stated by using (\ref{lim}).  For any fixed measurable function $\a(t):[0,\infty)\to A$ (control process), 
let  
 $(X_t^{\a}, A_t^{\a})$ be 
the stochastic process defined by 
\begin{eqnarray}\label{process} 
X_{t}^{\a}&=&x+\int_0^t \sigma^{\a} (X_s^{\a})dW_s+
\int_0^t b^{\a}(X_s^{\a}) ds-\int_0^t \g(X_s^{\a})dA_s \qquad t\ge 0 \nonumber,\\ 
A_t^{\a}&=&\int_0^t 1_{\p\Omega}(X_s^{\a})dA_s \quad \hbox{is continuous, non decreasing in}\quad t\ge 0, 
\end{eqnarray}
where $b^{\a}=(b_i^{\a})_i$, $1_{\p\Omega}(\cdot)$ a characteristic function on $\p\Omega$, 
$W_t$ $(t\ge 0)$ an $m-$dimensional Brownian motion. The study of the existence and  the uniqueness of $(X_t^{\a},A_t^{\a})$ is called the Skorokhod problem, and its solvability   is known under the preceding assumptions. We refer the readers to P.-L. Lions and A.S. Sznitman \cite{lsz}, P.-L. Lions, J.M. Menaldi and A.S. Sznitman \cite{lpms}, and P.-L. Lions \cite{new}.
  Let
$$
	J_{\l}^{\a} (x)=E_x \int_0^{\infty} e^{-\l t}g(X_t^{\a})
	1_{\p\Omega}(X_t^{\a}) dA_t,
$$
and define 
\begin{equation}\label{ul}
u_{\l}(x)=\inf_{\a(\cdot)} J_{\l}^{\a}(x) \qquad \hbox{in}\quad \Omega,
\end{equation}
where the infimum is taken over all possible control processes. 
It is known that $u_{\l}$ 
is the unique solution of (\ref{le})-(\ref{lb}). (See, P.-L. Lions and N.S. Trudinger \cite{lt}, and M.I. Freidlin and A.D. Wentzell \cite{fw}.)
Thus,
\begin{equation}\label{long}
d=\lim_{\l\downarrow 0} \inf_{\a(\cdot)} \l E_x \int_0^{\infty} 
e^{-\l t}g(X_t^{\a}) 1_{\p\Omega}(X_t^{\a}) dA_t,
\end{equation}
 if the right hand side of (\ref{lim}) exists, which represents the fact that 
 the number $d$ is  the long time averaged 
reflection force on the boundary. (Each time the tragectory reaches to $\p\Omega$, it 
gains the force $g(x)$ and is pushed back in the direction of $-\g(x)$.) We remark the similarity 
 of the convergence (\ref{lim}) to the so-called ergodic problem for HJB equations. That is, by considering,   
$$
	\l u_{\l}(x)+F(x,\n u_{\l}, \H u_{\l})=0 \qquad \hbox{in} \quad \Omega,
$$ 
$$
	<\n u_{\l}(x), \g (x)>=0 \qquad \hbox{on} \quad \p\Omega,
$$
it is known that an unique number $d'$ exists such that
$$
	\lim_{\l \downarrow 0} \l u_{\l}(x)=d' \qquad \hbox{uniformly in}\quad \Omega. 
$$
We refer the readers to M. Arisawa and P.-L. Lions \cite{apl}, M. Arisawa \cite{ar1},  \cite{ar2}, A. Bensoussan \cite{ben}  for the various types (operators and boundary conditions) of ergodic problems. As the above ergodic problem "in the domain", the existence of $d$ in (\ref{2}) "on the boundary" relates to the ergodicity of the stochastic process (\ref{process}). 
Even for some classes of degenerate elliptic operators $F$, the number $d$ in (\ref{2}) exists. 
We remark this in Section 4, below.\\

	Next, we turn our interests to the homogenization. The unique 
existence of $d$ in (\ref{1})-(\ref{2}) plays an essential role to study  the homogenization of oscillating 
Neumann boundary conditions. The simplest example is as follows.

\begin{example}{\bf Example 1.2.} Let  $c$, $g$, $f_1$$(x,\xi_1)$ 
be functions defined in $(x,\xi_1)\in \R^2\times {\bf R}\backslash {\bf Z}$ (periodic in $\xi_1$ with period $1$). Assume that $f_1\ge 0$, and that there exists a constant $c_0>0$ such that  $c>c_0>0$. For any $\e\ge 0$, let
$$
\Omega_{\e}=\{(x_1,x_2)|
\quad \e f_1(x,\frac{x_1}{\e})\le x_2\le b,\quad |x_1|\le a\},
$$
$$
\G_{\e}=\{(x_1,x_2)|\quad x_2=\e f_1(x,\frac{x_1}{\e}) \} \cap \p\Omega_{\e}.
$$
Let $u_{\e}(x)$ $(\e>0)$ be the solution of 
\begin{equation}\label{lap}
-\Delta u_{\e}=0 \qquad \hbox{in} \qquad \Omega_{\e}, 
\end{equation}
\begin{equation}\label{be}
<\n u_{\e}(x),{\bf n}_{\e}(x)>+c(x,\frac{x_1}{\e})u_{\e}=g(x,\frac{x_1}{\e}) \qquad \hbox{on} 
 \quad \G_{\e},
\end{equation}
\begin{equation}\label{de}
u_{\e}=0 \qquad \hbox{on} \quad \p\Omega_{\e}\backslash \G_{\e},
\end{equation}
where ${\bf n}_{\e}(x)$ is the outward unit normal to $\G_{\e}$. Then, as $\e\downarrow 0$, $u_{\e}$ 
converges to a unique functiont $u(x)$ uniformly in $\overline{\Omega_0}$, which is the solution of 
$$
-\Delta u=0 \qquad \hbox{in} \quad \Omega_0,
$$
\begin{equation}\label{ol}
<\n u(x), \nu (x)> + \overline{L}(x,u,\n u)=0 \qquad \hbox{on} \quad \G_0,
\end{equation}
$$
u=0 \qquad \hbox{on} \quad \p\Omega_0 \backslash \G_0,
$$
where $\nu $ is the outward unit normal to $\G_0$, and $\overline{L}$ is defined as 
follows. \\
 Let $O(x)=\{(\xi_1,\xi_2)|\quad \xi_2\ge f_1(x,\xi_1), \quad \xi_1\in \R\backslash {\bf Z} \}$. 
Then, for any fixed $(x,r,p)$$\in \Omega\times \R\times \R^2$, there exists a unique number $d(x,r,p)$ such that 
$$
-\Delta_{\xi}v\equiv 
-(\frac{\p^2 v}{\p \xi_1^2}+\frac{\p^2 v}{\p \xi_2^2})=0 \qquad \hbox{in} \quad O(x),
$$
$$
d(x,r,p)+<\n_{\xi} v, \gamma(\xi)>-(\sqrt{1+(\frac{\p f_1}{\p \xi_1})^2}g-
\sqrt{1+(\frac{\p f_1}{\p \xi_1})^2}cr-p_1\frac{\p f_1}{\p \xi_1})=0
$$
 on $\p O(x)$, where $\gamma (\xi)=(\frac{\p f_1}{\p \xi_1},-1)$ ($\xi\in \p O(x)$),  and 
\begin{equation}\label{effe}
\overline{L}(x,r,p)=-d(x,r,p).
\end{equation}
\end{example}

In A. Friedman, B. Hu, and Y. Liu \cite{af}, a similar problem to the above example (linear, three scales case) was treated by the variational approach. 
 (See also \cite{che}.) We shall extend  the result
(including Example 1.2.) to nonlinear  problems by using the existence of the long time averaged reflection number $d$   
in (\ref{1})-(\ref{2}). As Example 1.2 indicates, the effective limit boundary condition (\ref{ol}) is defined by using the long time averaged number in $(\ref{effe})$.   Our present approach was inspired by the classical method of formal asymptotic expansions of A. Bensoussin, J.L. Lions, and G. Papanicolaou \cite{jl2}.  This approach is closely related to  the ergodic problem for HJB equations described in the preceding part of this introduction. For  the application of the ergodic problem ( \cite{apl}, \cite{ar1}, \cite{ar2})  to obtain the effective P.D.E. in the domain, we refer the 
readers to M. Arisawa \cite{ar3}, \cite{ar4}, M. Arisawa and Y. Giga \cite{ag}, L.C. Evans \cite{ev1}, \cite{ev2}, and 
P.-L. Lions, G, Papanicolaou, and S.R.S. Varadhan \cite{lpv}n. As far as we know, there is no existing reference which treats the homogenization of the oscillating Neumann 
boundary conditions from the view point  of the ergodic problem.\\

	The plan of this paper is the following. \\
$\S 1.$ Introduction.\\
$\S 2.$ Existence and uniqueness of the number $d$ in the case of the 
bounded domain.\\
$\S 3.$ Existence and uniqueness of the number $d$ in the case of the 
half space.\\
$\S 4.$ Some remarks on the degenerate elliptic operators case.\\
$\S 5.$ Homogenization of the oscillating Neumann boundary conditions.\\ 

	Throughout of this paper, the gradient and the Hesse matrix of $u(x)$ ($x\in \Omega\subset {\bf R}^n$) ( resp. $v(\xi)$ ($\xi\in \Omega'\subset {\bf R}^n  $)) are denoted by  $\n u(x)$, $\n^2 u(x)$ ( resp. $\n_{\xi} v(\xi)$, $\n^2_{\xi} v(\xi)$ or $D^2_{\xi} v(\xi)$). For $u(x)$ ($x\in \Omega\subset {\bf R}^n$), the partial derivatives  in $x_i$, $x_j$ ($1\le i,j\le n$) are denoted by $\frac{\p u}{\p x_i}=D_i u$, $\frac{\p^2 u}{\p x_i\p x_j}=D_{ij} u$, etc., and the derivatives in the directions of 
$y,z\in {\bf R}^n$ are denoted by $D_y u=\sum_{i=1}^n y_i\frac{\p u}{\p x_i}$, 
 $D_{yz} u=\sum_{i,j=1}^n y_i z_j \frac{\p^2 u}{\p x_i\p x_j}$, etc.. 
 When a function $w(x,\xi)$ depends on both variables of $x\in {\bf R}^n$ and $\xi\in {\bf R}^n$,  and when we consider the derivatives $\frac{\p^2 w(x,\xi)}{\p x_k \p \xi_l}$ etc., we denote them by $D_{ij}w(x,\xi)$ ($1\le i,j\le 2n$), etc.. 
For the twice continuously differentiable function $u(x)$ ($x\in \Omega\subset {\bf R}^n$), we 
denote $|u|_{L^{\infty}(\Omega)}$$=\sup_{x\in \Omega}|u|$, $|\n u|_{L^{\infty}(\Omega)}$$=\sup_{x\in \Omega}\sup_{1\le i\le n}|\frac{\p u}{\p x_i}(x)|$, 
$|\n^2 u|_{L^{\infty}(\Omega)}$$=\sup_{x\in \Omega}\sup_{1\le i,j\le n}|\frac{\p^2 u}{\p x_i\p x_j}(x)|$,
\begin{eqnarray}
	&&|u|_{\beta;\Omega}=\sup_{(x,y)\in \Omega\times\Omega}
	\frac{|u(x)-u(y)|}{|x-y|^{\beta}},\quad
	|\n u|_{\beta;\Omega}=\sup_{1\le i\le n} \sup_{(x,y)\in \Omega\times\Omega}
	\frac{|\frac{\p u}{\p x_i}(x)-\frac{\p u}{\p x_i}(y)|}{|x-y|^{\beta}}\quad 	0<\beta\le 1,\nonumber\\
	&&|u|_{j,\beta;\Omega}=|\n^j u|_{L^{\infty}(\Omega)}+
	 \sup_{x\neq y\in \Omega}\frac{|\n^j u(x)-\n^j u(y)|}{|x-y|^{\beta}}\qquad 		0<\beta\le 1, \quad j=1,2.\nonumber
\end{eqnarray}
 We consider the solvability of PDEs in the framework of viscosity solutions, and treat the second-order sub and super differentials of upper and lower 
semi continuous functions $u(x)$ and $v(x)$ ($x\in D\subset {\bf R}^n$) at a point $\overline{x}$ in the domain $D$. We denote them by $J^{2,+}_{D}u(\overline{x})$ (the second-order superjets of $u$ at $\overline{x}$) and 
$J^{2,-}_{D}v(\overline{x})$ (the second-order subjets of $v$ at $\overline{x}$) respectively. (See M.G. Crandall and P.-L. Lions \cite{cpl}, M.G. Crandall, H. Ishii and P.-L. Lions \cite{users}, and W.H. Fleming and H.M. Soner n\cite{flso}.) 
 We use the notation $B(x,r)$ $(x\in \Omega,\quad r>0)$ for the open ball centered at $x$ 
with radius $r>0$. \\

	The author is grateful to Professors P.-L. Lions, H. Ishii and Y. Giga for 
 their helpful comments and encouragements. She thanks to Professors O. Alvarez, G. Barles and T. Mikami for the  discussions and suggestions on this subject. Finally, she also thanks to anonymous referee for his kind comments and interesting suggestions.\\

\section{Existence and uniqueness of the long time 
 averaged reflection force in the bounded domain.}

	In this section, the existence and uniqueness of the number $d$ in (\ref{1})-(\ref{2}) is shown in the case that $\Omega$ satisfies (\ref{doma}). The 
Hamiltonian  
$F(x,\n u, \H u)$, given in (\ref{hjb}), positively homogeneous in degree one, is 
assumed to  satisfy 
(\ref{unif}) and (\ref{lips}); the vector field $\g$ on $\p\Omega$ is assumed to satisfy  (\ref{gamma}). For the existence, we further assume that 
\begin{equation}\label{bdd}
|a_{ij}^{\a}, |\n a_{ij}^{\a}|,|\H a_{ij}^{\a}|,|b_i^{\a}|, |\n b_i^{\a}|, |\n^2 b_i^{\a}|\le K
\qquad \hbox{any} \quad x\in \Omega, 
\quad 1\le  i,j\le n, \quad \a\in A, 
\end{equation}
where $K>0$ is a constant, and that 
$\g$, $g$ can be extendable in a neighborhood $U$ of $\p\Omega$ to twice continuously differentiable functions so that 
\begin{equation}\label{gbd}
|\n \g|,|\H \g|,|\H g|,|\H g|\le K \qquad \hbox{in}\quad U,
\end{equation}
where $K>0$ is the constant in (\ref{bdd}). For the existence of $d$, we approximate (\ref{1})-(\ref{2}) 
by (\ref{le})-(\ref{lb}) ($\l\in (0,1)$) and examine the 
regularity of $u_{\l}$, uniformly in $\l$.  In order to have (\ref{lim})-(\ref{ulim}), we need the following estimates.

\begin{proposition}{\bf Theorem 2.1.} Assume that $\Omega$ is (\ref{doma}), and that 
(\ref{unif}), (\ref{gamma}), (\ref{bdd}) and (\ref{gbd}) hold. 
 Then there exists a unique solution $u_{\l}\in C^{1,1}(\overline{\Omega})\cap C^{2,\beta}(\Omega)$ of (\ref{le})-(\ref{lb}), where $\beta>0$ depends on $n$ and $\Lambda_1/\l_1$. Moreover for any  fixed $x_0\in \Omega$, there exists a constant $C>0$ such that the following estimates hold.
\begin{equation}\label{ub}
|u_{\l}-u_{\l}(x_0)|_{L^{\infty}(\overline{\Omega})}
\le C \qquad \hbox{any} \quad  \l \in (0,1),
\end{equation} 
\begin{equation}\label{umb}
|\n u_{\l}|_{L^{\infty}(\overline{\Omega})}\le C \qquad \hbox{any}\quad  \l \in (0,1),
\end{equation}
\begin{equation}\label{c11}
|\n u_{\l}|_{1;\overline{\Omega}}\le C \qquad \hbox{any}\quad  \l \in (0,1).
\end{equation}
\end{proposition}

{\bf Remark 2.1} One can replace the conditions (\ref{bdd})-(\ref{gbd}) to other conditions 
 ,for example those in \cite{ipl}n, to have 
$$
	|u_{\lambda}(x)-u_{\lambda}(y)|\le C|x-y|^{\theta}\qquad \hbox{any}\quad x,y\in \overline{\Omega},\quad \lambda\in (0,1),
$$
where $C>0$, $\theta\in (0,1)$ are independent on $\lambda>0$. The proof of this inequality 
 can be done in a similar way to \cite{ipl}, but by taking account of the Neumann type boundary 
conditions, and also by using the estimate (\ref{ub}). We do not write the proof in this direction here, but shall use the method in a future occassion. \\

{\it Proof of Theorem  2.1.} For each $\l>0$, the existence and uniqueness of $u_{\l}\in C^{1,1}(\overline{\Omega})\cap C^{2,\beta}(\Omega)$ is established in P.-L. Lions and N.S. Trudinger \cite{lt}n. We are to show the uniform  (in $\l\in (0,1)$) regularity (\ref{ub})-(\ref{c11}) in the following two steps. In Step 1, (\ref{ub}) will be shown, 
 and in Step 2, (\ref{umb}) and (\ref{c11}) will be  shown. \\ 

\underline{Step 1.} We prove (\ref{ub}) by a contradiction argument. 
Let $x_0\in \Omega$ be fixed. Assume, as $\l>0$ goes to $0$
$$
|u_{\l}-u_{\l}(x_0)|_{L^{\infty}(\overline{\Omega})}\to \infty.
$$
 Set  
$$
\e_{\l}\equiv |u_{\l}-u_{\l}(x_0)|_{L^{\infty}(\overline{\Omega})}^{-1} \qquad \l\in (0,1),
$$
and let $v_{\l}\equiv \e_{\l}(u_{\l}-u_{\l}(x_0))$. Then, 
$$
|v_{\l}|_{L^{\infty}(\overline{\Omega})}=1, \quad v_{\l}(x_0)=0 
\qquad \hbox{any} \quad  \l\in (0,1).
$$
From (\ref{hjb}), 
  $v_{\l}$ satisfies $F(x,\n v_{\l}, \H v_{\l})=0$ in $\Omega$, and from 
(\ref{unif})  the Krylov-Safonov inequality (see \cite{cabre}n for instance) leads: 
for any compact set $V\subset \subset \Omega$, there exists a constant $M_V>0$ such that 
\begin{equation}\label{int}
|\n v_{\l}|_{L^{\infty}(\overline{V})}\le M_V \qquad \hbox{any}\quad  \l \in (0,1).
\end{equation}
We denote 
$$
v^{\ast}(x)=\limsup_{\l\downarrow 0, y\to x}v_{\l}(y),
\quad v_{\ast}(x)=\liminf_{\l\downarrow 0, y\to x}v_{\l}(y).
$$
Then, since 
$v_{\lambda}(x_0)=0$ $(\forall \lambda\in (0,1))$, from (\ref{int}) we have 
\begin{equation}\label{con}
v^{\ast}(x_0)=v_{\ast}(x_0)=0,
\end{equation}
\begin{equation}\label{or}
|v^{\ast}|_{L^{\infty}(\overline{\Omega})}=1,\quad \hbox{or}\quad 
|v_{\ast}|_{L^{\infty}(\overline{\Omega})}=1.
\end{equation}
 From (\ref{2}), $v_{\l}$ satisfies 
$$
<\n v_{\l},\g(x)>=\e_{\l}g-\l(v_{\l}+\e_{\l}u_{\l}(x_0)),
$$
and by the comparison result for (\ref{le})-(\ref{lb})
$$
|\l u_{\l}(x_0)|_{L^{\infty}(\overline{\Omega})}\le C \qquad \hbox{any}\quad \l\in (0,1), 
$$
where $C>0$ is a constant.
 By letting $\l\downarrow 0$, $v^{\ast}$ and $v_{\ast}$ are viscosity solutions of 
\begin{equation}\label{sb}
<\n v^{\ast},\g(x)>\le 0 \qquad \hbox{on}\quad \p\Omega,
\end{equation}
\begin{equation}\label{sp}
<\n v_{\ast},\g(x)>\ge 0 \qquad \hbox{on}\quad \p\Omega,
\end{equation}
and $v^{\ast}(x)$ (resp. $v_{\ast}(x)$) ($x\in \Omega$) satisfies 
$$
\qquad   F(x,\n v^{\ast},\H v^{\ast})\le 0, \qquad(\hbox{resp.}\quad  F(x,\n v_{\ast},\H v_{\ast})\ge 0) \qquad \hbox{in}\quad \Omega. \qquad  \quad \hbox{(10)'}
$$ 
(We refer the readers to \cite{users} and G. Barles and B. Perthame \cite{bp}n for this 
stability result.)

	Now we employ  the strong maximum principle of  M. Bardi and F. Da-Lio \cite{bardi}.  Remark that $F(x,p,R)$ given in (\ref{hjb}), satisfying (\ref{unif})  and (\ref{bdd}) enjoys the following two properties of  (\ref{bardi1}) and (\ref{bardi2}). \\

	(Scaling property) For any  $x_0\in \Omega$, for any $\eta>0$, there exists 
a function $\phi$$:(0,1)\to(0,\infty)$ such that  
\begin{equation}\label{bardi1}  
\overline{F}(x,\xi p,\xi R)\ge \phi(\xi)\overline{F}(x,p,R) \qquad \hbox{any}\quad \xi\in (0,1),
\end{equation}
holds for any $x\in B(x_0,\eta)$, $0<|p|\le \eta$, $|R|\le \eta$.\\

	(Nondegeneracy property) For any $x_0\in \Omega$, for any small vector $\nu\neq 0$, there exists 
a positive number $r_0$ such that 
\begin{equation}\label{bardi2}
\overline{F}(x_0,\nu,I-r \nu\otimes\nu)>0 \qquad \hbox{any}\quad r>r_0.
\end{equation}

	We cite the following result for our present and later purposes.\\

\begin{lemma}{\bf Lemma A. (n\cite{bardi}n)} (Strong maximum priciple) Let 
$\Omega \subset {\bf R}^n$ be an open set and let $u$ be an upper semicontinuous viscosity subsolution of 
$$
\overline{F}(x,\n u,\H u)=0 \qquad \hbox{in}\quad \Omega, 
$$
which attains a maximum in $\Omega$. Assume that $\overline{F}$ 
satisfies (\ref{bardi1}), (\ref{bardi2}), and \\

 for any $x_0\in \Omega$ there exists $\rho_0>0$ such that for any $\nu\in 
B(0,\rho_0)\backslash \{0\}$, (\ref{bardi2})  
\begin{equation}\label{bardi3} \hbox{holds for some}\quad r_0>0.\qquad \qquad \qquad \qquad \qquad \qquad \qquad \qquad \qquad \qquad \qquad \qquad \quad 
\end{equation}  

	Then, $u$ is a constant. 
\end{lemma}

\bigskip

 We go back to the proof of (\ref{ub}). Assume that $|v^{\ast}|_{L^{\infty}(\overline{\Omega})}=1$ holds in (\ref{or}). (The another case
 of  $|v_{\ast}|_{L^{\infty}(\overline{\Omega})}=1$ can be treated similarly.) Thus from  (\ref{con}), $v^{\ast}$ is not constant, and from (10)' and the strong maximum principle (Lemma A), $v^{\ast}$ attains its maximum at a point $x_1\in \p\Omega$: 
$$
v^{\ast}(x_1)>v^{\ast}(x) \qquad \hbox{any}\quad x\in \Omega.
$$ 
 Since $\p\Omega$ is $C^{3,1}$, the interior sphere condition (see D. Gilbarg and N.S. Trudinger \cite{trudinger}n) is satisfied : there exists $y\in \Omega$ such that for $R=|x_1-y|$ 
$$
B(y,R)\in \Omega,\qquad x_1\in \p B(y,R). 
$$ 
Let 
$$
\phi(x)=e^{-cR^2}-e^{-c|x-y|^2} \qquad x\in \Omega,
$$
where $c>0$ is a constant large enough so that 
\begin{eqnarray}
	&&F(x_1,\n \phi(x_1),\n^2 \phi (x_1))\nonumber\\
	&&=
	F(x_1,2c(x_1-y)e^{-c|x_1-y|^2},2ce^{-c|x_1-y|^2}(I-2c(x_1-y)\otimes(x_1-y)))	\nonumber\\
	&&=2ce^{-c|x_1-y|^2}F(x_1,x_1-y,I-2c(x_1-y)\otimes(x_1-y))>0\nonumber
\end{eqnarray}
holds. (Here, we used (\ref{hjb}), (\ref{bardi2}) and (\ref{bardi3}).)
 By the lower semicontinuity of $F$ in $x$, there exists $r\in B(0,R)$ and $C'>0$ such that 
\begin{equation}\label{ish}
	F(x,\n \phi(x),\n^2 \phi(x))\ge C'>0 \qquad \hbox{in} \quad 
	B(x_1,r)\cap\overline{\Omega}. 
\end{equation} 
We claim that 
\begin{equation}\label{rel}
v^{\ast}(x)-v^{\ast}(x_1)-\phi(x)\le 0\qquad \hbox{in} \quad B(x_1,r)\cap 
\overline{\Omega}.
\end{equation}
 In fact, if $x\in B(y,R)^c$,  $\phi(x)\ge0$ and (\ref{rel}) holds. Assume that for 
 $x'\in B(x_1,r)\cap B(y,R)$ (\ref{rel}) does not hold, and 
$$ 
v^{\ast}(x')-v^{\ast}(x_1)-\phi(x')=\max_{B(x_1,r)\cap B(y,R)} v^{\ast}(x)-v^{\ast}(x_1)-\phi(x). 
$$
 Then by the definition of the viscosity solution,
$$
	F(x',\n \phi (x'),\n^2 \phi (x'))\le 0,
$$
which contradicts to (\ref{ish}). Therefore, (\ref{rel}) holds. By remarking that 
$\phi(x_1)=0$, (\ref{rel}) indicates that $v^{\ast}-\phi$ takes its maximum at $x_1\in \p\Omega$. Since $v^{\ast}$ satisfies (\ref{sb}) in the sense of viscosity solutions, 
either 
$$
<\phi(x_1),\g(x_1)>\le 0,
$$ 
or
$$
F(x_1,\n \phi(x_1),\n^2 \phi(x_1))\le 0
$$
must be satisfied. However from the definition of $\phi$, (\ref{gamma}) and (\ref{ish}), 
both of the above are not satisfied. We got a contradiction, and proved (\ref{ub}). \\

\underline{Step 2.} To obtain (\ref{umb}) and (\ref{c11}), we appply (\ref{ub}) in the argument of \cite{lt}n. First, we regularlize the Hamiltonian $F$.  Let $\rho$ be a mollifier on $\Rn$ ($\rho\ge 0$, $\rho \in C_0^{\infty}(\Rn)$, $\int \rho =1$). For any $\delta>0$, set
$$
	h_{\d}(y)=\d^{-n}\int_{\R^N}\rho(\frac{y-z}{\d})(\inf_{1\le k\le N}z_k) dz,
$$ 
$$
F_{\d}^N [u]\equiv h_{\d}(L^{\a_1}u,...,L^{\a_N}u),
$$
 where
$$
 L^{\a_l}u=-\sum_{i,j=1}^n a_{ij}^{\a_l}\frac{\p^2 u}{\p x_i\p x_j}
-\sum_{i}^{\a}b_{i}^{\a_l}\frac{\p u}{\p x_i} \qquad 1\le l\le N. 
$$
 Remark that for any $\d\in (0,1)$, the operator $F_{\d}^N(x,p,R)$ satisfies 
\begin{equation}\label{st1}
\l_1 I\le (\frac{\p F_{\d}^N}{\p r_{ij}}(x,p,R))_{1\le i,j\le n}\le \L_1 I \qquad x\in \Omega,\quad R\in {\bf S}^n,
\end{equation} 
\begin{equation}\label{st2}
F_{\d}^N(x,p,R)\le \mu_0(1+|p|+|R|) \qquad x\in \Omega,\quad R\in {\bf S}^n,
\end{equation}
\begin{equation}\label{st3}
|\frac{\p F_{\d}^N}{\p x}|, |\frac{\p F_{\d}^N}{\p p}|, |\frac{\p F_{\d}^N}{\p R}|
\le \mu_1\{(1+|p|+|R|)|x|+|p|+|R|\} \qquad x\in \Omega,\quad R\in {\bf S}^n,
\end{equation}
\begin{equation}\label{st4}
|\frac{\p^2 F_{\d}^N}{\p x^2}|, |\frac{\p^2 F_{\d}^N}{\p x\p p}|, 
|\frac{\p^2 F_{\d}^N}{\p x\p R}|\le \mu_2\{(1+|p|+|R|)|x|+|p|+|R|\}\times|x|\qquad x\in \Omega, \quad R\in {\bf S}^n,
\end{equation}
where $\mu_i$ ($i=0,1,2$) are positive constants, and $|p|=\max_{1\le i \le n} |p_i|$ 
 ($p=(p_i)_{1\le i \le n}$), $|R|=\max_{1\le i,j\le n}|r_{ij}|$ ($R=(r_{ij})_{1\le i,j\le n}$). \\

We need the following a priori estimates.

\begin{theorem}{\bf Lemma 2.2.} 
	Let $u_{\l,N}^{\d}\in C^4(\Omega)\cap C^3(\overline{\Omega})$ be 
a solution of 
\begin{equation}\label{three}
	F_{\d}^N(x,\n u_{\l,N}^{\d}, \H u_{\l,N}^{\d})=0 \qquad \hbox{in}\quad \Omega,
\end{equation}
\begin{equation}\label{bthree}
	\l u_{\l,N}^{\d}+<\n u_{\l,N}^{\d},\g(x)>-g(x)=0\qquad \hbox{on}\quad \p\Omega.
\end{equation}
 Then, there exists $C>0$ such that  
\begin{equation}\label{tha}
	|\n u_{\l,N}^{\d}|_{L^{\infty}(\overline{\Omega})},\quad 
	|\n^2 u_{\l,N}^{\d}|_{L^{\infty}(\overline{\Omega})}\le C 
	\qquad \hbox{any} \quad \delta, \quad \l\in (0,1), \quad N\in{\bf N},
\end{equation}
where $C>0$ depends on $n$, $\l_1$, $\L_1$, $\mu_i$ ($i=0,1,2$),  $\Omega$, and  $K$. 
\end{theorem}

{\bf Remark 2.2.} In the estimates of \cite{lt}, Theorem 2.1n,  the above constant $C$ 
 depends also on $\l\in(0,1)$. \\

\bigskip
	By delaying the proof of Lemma 2.2, we shall show how (\ref{tha}) leads (\ref{umb}) and (\ref{c11}). By the method of continuity, for each $\d>0$	
the a priori estimate (\ref{tha}) yields the existence of $u_{\l,N}^{\d}$ 
$\in C^3(\Omega)$$\cap C^{2,\a}(\overline{\Omega})$ of (\ref{three})-(\ref{bthree}).
 Put $w_{\l,N}^{\d}=u_{\l,N}^{\d}-u_{\l,N}^{\d}(x_0)$. The same argument as in Step 1 works for $w_{\l,N}^{\d}$, and 
$$
	|w_{\l,N}^{\d}|_{L^{\infty}(\overline{\Omega})}\le C 
	\qquad \hbox{any}\quad \d,\quad \l\in (0,1), \quad N\in {\bf N}.
$$ From (\ref{tha}), by extracting a subsequence of $\d'\downarrow 0$, there exists $w_{\l,N}$ $\in C^{1,1}(\overline{\Omega})$ such that 
$$
	\lim_{\d'\downarrow 0} w_{\l,N}^{\d}=w_{\l,N} \qquad \hbox{uniformly in}\quad 
	\overline{\Omega},
$$
$$
	\lim_{\d'\downarrow 0} \n w_{\l,N}^{\d}=\n w_{\l,N} \qquad \hbox{uniformly 		in}\quad \overline{\Omega},
$$
and 
$$
	|w_{\l,N}|_{L^{\infty}(\overline{\Omega})}, \quad 	|\n w_{\l,N}|_{L^{\infty}(\overline{\Omega})}, \quad 		|\n w_{\l,N}|_{1;\overline{\Omega}}\le C \qquad 
	\hbox{any} \quad \l\in (0,1),\quad N>0.  
$$
	On the other hand, from (\ref{st1}) and  the Evans-Krylov interior 
estimate (see, e.g. L.C. Evans \cite{ev0}, X. Cabre and L.A. Caffarelli \cite{cabre},   N.V. Krylov \cite{kr1}, \cite{kr2}, and  \cite{lt}n,) leads 
for any $\Omega'\subset \subset \Omega$ 
$$
	|\n^2 w_{\l,N}^{\d}|_{\a;\Omega'}\le C\qquad \hbox{any}\quad \d\in (0,1),
$$
where $C>0$ depends on $\Omega'$ and $\a\in (0,1)$. Thus, we obtain $w_{\l,N}\in$$C^{1,1}(\overline{\Omega})\cap C^{2,\beta}(\Omega)$ of 
$$
	\max_{1\le l\le N}\{L^{\a_l}w_{\l,N}\}=0\qquad \hbox{in}\quad \Omega,
$$
$$
	\l w_{\l,N}+<\n w_{\l,N},\g(x)>-g(x)=0\qquad \hbox{on}\quad \p\Omega.
$$
	Letting $N\to \infty$, we obtain (\ref{umb}) and (\ref{c11}) from the preceding 
estimates. \\

	In the following, we shall prove Lemma 2.2.\\

{\it Proof of Lemma 2.2.}
 Set 
\begin{equation}\label{vm}
v_{\l,N}^{\d}\equiv \frac{u_{\l,N}^{\d}-u_{\l,N}^{\d}(x_0)}{|\n (u_{\l,N}^{\d}-u_{\l,N}^{\d}(x_0)) |_{L^{\infty}(\overline{\Omega})}}.
\end{equation}
From (\ref{ub}), there exists a constant $M_1>0$ such that 
\begin{equation}\label{vM}
|v_{\l,N}^{\d}|_{L^{\infty}(\overline{\Omega})},\quad |\n v_{\l,N}^{\d}|_{L^{\infty}(\overline{\Omega})}\le M_1 \qquad \hbox{any} \quad 
\d,\quad \l\in (0,1), \quad N\in {\bf N}.
\end{equation}
It is clear that
\begin{equation}\label{ve}
F_{\d}^N(x,\n v_{\l,N}^{\d},\H v_{\l,N}^{\d})=0 \qquad \hbox{in}\quad \Omega,
\end{equation}
\begin{equation}\label{vb}
\l v_{\l,N}^{\d}+<\n v_{\l,N}^{\d},\g(x)>-\overline{g}=0 \qquad \hbox{on} \quad
 \p\Omega,
\end{equation}
where 
$$
\overline{g}=\frac{g-\l u_{\l,N}^{\d}(x_0)}{|\n (u_{\l,N}^{\d}-u_{\l,N}^{\d}(x_0))|_{L^{\infty}(\overline{\Omega})}}. 
$$

We need the following Lemma. 

\begin{lemma}{\bf Lemma 2.3.}
Let $v_{\l,N}^{\d}$ be defined in (\ref{vm}). Then, there exists $C>0$ such that
\begin{equation}\label{vd}
	|\H v_{\l,N}^{\d}|_{L^{\infty}(\overline{\Omega})}\le C \qquad 
	\hbox{any}\quad \d,\quad \l\in (0,1),\quad N\in {\bf N}.
\end{equation}
\end{lemma}

	Lemma 2.3 will lead our present goal (\ref{tha}) in Lemma 2.2.  In fact,  from (\ref{vm}), (\ref{vd}), we have
\begin{equation}\label{gs}
\sup_{\overline{\Omega}}|\H u_{\l,N}^{\d}|\le C(1+\sup_{\overline{\Omega}}|\n u_{\l,N}^{\d}|).
\end{equation}
 We use the following interpolation inequality in the above. 

\begin{lemma}{\bf Lemma B. (\cite{trudinger}, Lemma 6.35)} Suppose $j+\beta<k+\a$, where $j=0,1,2,...$; $k=1,2,...$, and $0\le \a,\beta\le 1$. Let $D$ be a $C^{k,\a}$ domain 
 in $\Rn$, and assume $u\in C^{k,\a}(\overline{D})$. Then, for any $\e>0$ and some 
constant $C=C(\e,j,k,D)$ we have 
$$
|u|_{j,\beta;D}\le C|u|_{L^{\infty}(D)}+\e|u|_{k,\a;D}.
$$ 
\end{lemma}
  
	By putting $j=1$, $k=2$, $\a=\beta=0$ in Lemma B,  (\ref{gs}) leads (\ref{tha}) in 
Lemma 2.2. Finally,  we are to prove Lemma 2.3.\\

{\it Proof of Lemma 2.3.} For simplicity,  write $F=F_{\d}$, 
$v=v_{\l,N}^{\d}$. First, we examine the regularity of $v$ on $\p\Omega$. 
 By differentiating (\ref{ve}) twice with respect to a vector $\xi\in \Rn$, $|\xi|=1$, 
$$
\sum_{i,j=1}^n \frac{\p F}{\p r_{ij}}\frac{\p^2}{\p x_i \p x_j}D_{\xi} v + \sum_{i=1}^n \frac{\p F}{\p p_i}\frac{\p}{\p x_i}D_{\xi}v + \frac{\p F}{\p \xi}=0,
$$
$$
\sum_{i,j=1}^n \frac{\p F}{\p r_{ij}}\frac{\p^2}{\p x_i \p x_j}D_{\xi \xi} v + \sum_{i=1}^n\frac{\p F}{\p p_i}\frac{\p}{\p x_i}D_{\xi \xi}v + F_{\overline{X} \overline{X}}=0,
$$
where $F_{\overline{X} \overline{X}}$ is the derivarives of $F$ with respect to 
$\overline{X}=(\xi, \n (D_{\xi}v), \n^2 (D_{\xi}v))$. Using the structure conditions 
(\ref{st1})-(\ref{st4}), we obtain from above inequalities
\begin{equation}\label{A1}
|\sum_{i,j=1}^n \frac{\p F}{\p r_{ij}}\frac{\p^2}{\p x_i \p x_j}D_{\xi}v|\le C(1+|\H v|),
\end{equation}
\begin{equation}\label{A2}
\sum_{i,j=1}^n \frac{\p F}{\p r_{ij}}\frac{\p^2}{\p x_i \p x_j} D_{\xi \xi}v
\le C(1+|\H v|+|\H D_{\xi}v|),
\end{equation}
where $C>0$ depends on $n$, $M_1$, $\mu_1$ and $\mu_2$.  
	By the usual 
argument of flattening the boundary, we may assume that $\p\Omega$$=\{(x',x_n)|x_n\ge 0\}$ in a neighborhood of $x=0 \in \p\Omega$. Although by the change of variables, (\ref{ve})-(\ref{vb}) is transformed into $\overline{F}=0$ ($\overline{F}$ is the new Hamiltonian) etc., we keep to denote 
$\overline{F}=F$, etc., for simplicity.  
Denote $B_{r}^{+}=\{x\in B(0,r)|x_n>0\}$, 
 and for $\xi=(\xi_1,...,\xi_{n-1},0)\in{\bf R}^{n-1}$, $|\xi|\le 1$, consider 
\begin{equation}\label{w}
w(x,\xi)\equiv \eta^2(x,\xi)(z(x,\xi)+Av'),
\end{equation}
where $\eta$ is a smooth cut-off function to be precised in below, $A$ a constant, 
$$
z(x,\xi)\equiv D_{\xi\xi}v(x)=\sum_{ij}\frac{\p^2 v}{\p x_i\p x_j}\xi_i\xi_j,
\qquad v'\equiv \sum_{i=1}^{n-1} |\frac{\p v}{\p x_i}|^2. 
$$ 
By introducing  (\ref{st1}), (\ref{st2}), (\ref{vM}) and (\ref{ve}) into (\ref{A1}), we obtain 
$$
\sum_{i,j=1}^{n} (\frac{\p F}{\p r_{ij}}\frac{\p^2 z}{\p x_i\p x_j}+C_{ij}\frac{\p^2}{\p x_i \p x_j} D_{\xi}v)\le 
 C(1+|\H v|')
$$
where the coefficients $C_{ij}$ are such that $C_{in}=0$, $|C_{ij}|\le C$ 
 depending  on $n$, $\l_1$, $\mu_i$ ($i=0,1,2$), $M_1$, and 
$|\H v|'=(\sum_{i+j<2n}|\frac{\p^2 v}{\p x_i \p x_j}|^2)^{\frac{1}{2}}$. 
 Using the relations 
$$
\frac{\p}{\p x_i} D_{\xi_j}z=2\frac{\p^2}{\p x_i \p x_j} D_{\xi}v,
\quad D_{\xi_i \xi_j}z=2\frac{\p^2 v}{\p x_i \p x_j},
$$
we can take constants $C_0$ and $C$ such that the following 
 $(2n-1)\times(2n-1)$ matrix $(F_{ij}')_{ij}$:
$$
\sum_{i,j=1}^{2n-1}F_{ij}'D_{ij}z\equiv \sum_{i,j=1}^{n}\frac{\p F}{\p r_{ij}}\frac{\p^2 z}{\p x_i\p x_j}+
\frac{1}{2}\sum_{i=1}^{n}\sum_{j=1}^{n-1}C_{ij}\frac{\p}{\p x_i}D_{\xi_j}z
+C_0\sum_{j=1}^{n-1}D_{\xi_j \xi_j} z
$$
$$
	\le C(1+|\H v|')
$$
 is uniformly elliptic with minimum eigenvalue $\l'\ge \frac{\l_1}{2}$. 
 From (\ref{A1}), 
$$
\sum_{k=1}^{n-1} \frac{\p F}{\p r_{ij}}
\frac{\p^2 v}{\p x_i \p x_k}\frac{\p^2 v}{\p x_j \p x_k}
+\frac{1}{2}\frac{\p F}{\p r_{ij}}
\frac{\p^2 v'}{\p x_i \p x_j}\le C(1+|\H v|').
$$
 By combining the above two inequalities, we arrive at 
\begin{equation}\label{ineq}
\eta^2 \sum_{i,j=1}^{2n-1} F_{ij}' D_{ij}w
-2\sum_{i,j=1}^{2n-1} F_{ij}' D_{i}\eta^2 D_j w
\le
\end{equation}
$$
-2K\l (|\H v|')^2\eta^4+6(\sum_{i,j=1}^{2n-1} F_{ij}' D_{i}\eta D_j \eta)w
-2\eta(\sum_{i,j=1}^{2n-1} F_{ij}' D_{ij} \eta)w
-C(1+K)\eta^4(1+|\H v|')
$$
$$
\le -A\l w^2+C_A,
$$
where the constant $C_A$ depends on $n$, $\l_1$, $\mu_i$ ($i=0,1,2$) and $M_1$. 
(Remark that $C_A$ does not depend on $\l\in (0,1)$, for we have not yet used the boundary condition (\ref{vb})). \\
	Next, by differentiating (\ref{vb}) in the direction of $\xi_k$, $\xi_l$, 
\begin{equation}\label{bv}
\l D_{\xi_k}v+<\n (D_{\xi_k}v),\g>+<\n v,D_{\xi_k}\g>=D_{\xi_k}\overline{g}, 
\end{equation}
\begin{equation}\label{bbv}
\l D_{\xi_k \xi_l}v+<\n (D_{\xi_k \xi_l}v),\g>+<\n (D_{\xi_k}v),D_{\xi_l}\g>+<\n (D_{\xi_l}v),D_{\xi_k}\g> \qquad 
\end{equation}
$$
\qquad \qquad \qquad \qquad +<\n v,D_{\xi_k \xi_l}\g>=D_{\xi_k \xi_l}\overline{g}.
$$
Since 
$$
\frac{\p w}{\p x_i}=2\frac{w}{\eta}\frac{\p \eta}{\p x_i}+\eta^2(\frac{\p z}{\p x_i}+A\frac{\p v'}{\p x_i}),
$$
$$
\l w+<\n w,\g>-2\frac{w}{\eta} <\n \eta,\g>\qquad \qquad \qquad \qquad \qquad \qquad \qquad 
$$
$$
=\eta^2<\n z,\g>+\eta^2 A<\n v',\g>+\l \eta^2(z+Av'),
$$
and from (\ref{bbv}), 
$$
=\eta^2A<\n v',\g>+\l\eta^2 Av'-\eta^2<\n v,D_{\xi_k \xi_l}\g>-2 \eta^2
<\n (D_{\xi_k}v),D_{\xi_k}\g>. 
$$
From (\ref{gbd}) and (\ref{vM}), 
$$
|v'|,\quad |D_{\xi_k} \g|, \quad |D_{\xi_k \xi_l} \g|\le K \qquad \hbox{any} \quad 
1\le k,l\le n-1,
$$
and by (\ref{bv}) $<\n v',\g>$ and $<\n (D_{\xi_k}v),D_{\xi_k}\g>$ are bounded. 
 Therefore, we can fix $A$ so that 
$$
\l w+<\n w,\g>-2\frac{w}{\eta} <\n \eta,\g>\le C_1 \eta^2,
$$
where $C_1>0$
 depends on $n$, $\l_1$, $\mu_i$ ($i=0,1,2$), $K$ and 
$M_1$. (In particular, $C_1$ is independent of $\l \in (0,1)$.) 
 Now, fix
$$
	\eta(x,\xi)=[1-4\{|x'|^2+(x_n-\overline{\e}r)^2\}/r^2-|\xi|^2\}^{+},
$$
where for
$$
T=\{x\in B_r, \quad x_n=0\}, \qquad N=\{(x,\xi)\in \R^{2n-1}|\quad \eta(x,\xi)>0\},
$$
$$\overline{\e}=\zeta/\sqrt{1+\zeta^2}, \qquad \zeta=\sup_{T}\frac{|\g|}{\g_n}\le C. 
$$ 
 Then, on $T\cap \p N\cap \{w\ge 0\}$
$$
	<\n w,\g>+\l w\le C_2,
$$ 
where $C_2$ is independent of $\l \in (0,1)$. We take $\overline{w}=w+C_3{\l_1}^{-1}{x_n}$ 
so that 
$$
<\n \overline{w},\g>=<\n w,\g>+\g_n\frac{C_3}{\l_1}\g\le C_2-\l w+\g_n\frac{C_3}{\l_1}\le 0.
$$
 From the definition of $w$, the above constant $C_3$ can be taken uniformly in 
$\l\in(0,1)$. By applying the maximum principle to $\overline{w}$, instead of $w$, we obtain 
\begin{equation}\label{xx}
D_{\xi\xi}v(0)\le C,
\end{equation}
for any $\xi=(\xi_1,...,\xi_{n-1},0)$ ($|\xi|=1$), where $C>0$ depends only on 
$\eta$, $\l_1$, $\mu_i$ ($i=0,1,2$), $M_1$, $\Omega$ and $K$. ( $C$ is independent of  $\l\in(0,1)$.) As for the remaining inequalities, the same argument in \cite{lt}n is available. That is, by regarding 
$$
G(x)=\l v+<\n v,\g>-g(x)
$$ as a function in $B(0,r)$ ($0\in \p\Omega$, $\g$ and $g$ are extendable to some    neighborhood of $\p\Omega$ (\ref{gbd})), 
$$
|\sum_{i,j=1}^n \frac{\p F}{\p r_{ij}}\frac{\p^2 G}{\p x_i\p x_j}|\le C(1+M_2) \qquad  (M_2=\sup_{\Omega}|\H v|) \quad 
\hbox{in}\quad B(0,r),
$$
$$
G=0 \qquad \hbox{on}\quad \p\Omega,
$$
where $C$ depends on $n$, $M_1$, $\mu_1$, $K$, and does not depend on $\l\in(0,1)$. 
From this, the barrier argument leads 
\begin{equation}\label{bar}
|DG(0)|\le C\sqrt{1+M_2},
\end{equation} 
and we can extend the inequality (\ref{xx}) to 
\begin{equation}\label{rest}
D_{\xi\xi}v(0)\le C \qquad \hbox{any}\quad  |\xi|=1,\quad \xi\in{\bf R}^n.
\end{equation}
\newpage

 Then, by  the uniform ellipticity (\ref{st1}), the usual argument leads 
\begin{equation}\label{ryo}
\sup_{\p\Omega}|\H v|\le C\qquad \hbox{any}\quad  |\xi|=1,\quad \xi\in{\bf R}^n,
\end{equation}
where $C$ is independent of $\l\in(0,1)$. From (\ref{st1}), by coupling (\ref{ryo}) with the global Dirichlet bound for 
(\ref{ve})-(\ref{vb}) leads (\ref{vd}), and Lemma 2.3 was proved.\\

We complete the  proof of Theorem 2.1.\\

\begin{theorem}{\bf Theorem 2.4.} Assume that $\Omega$ is (\ref{doma}), and that (\ref{unif}), (\ref{gamma}), (\ref{bdd}) and (\ref{gbd}) hold. Then there exists a number $d$ and a function $u(x)\in C^{1,1}(\overline{\Omega})\cap C^{2,\alpha}(\Omega)$ ($\a\in (0,1)$) which satisfy (\ref{1})-(\ref{2}). 
\end{theorem}

{\it Proof of Theorem 2.4.}	From (\ref{ub})-(\ref{c11}) and the Evans-Krylov estimate, we can extract a subsequence $\l'\downarrow 0$ such that there exist a number $d$ and $u(x)\in$ 
$C^{1,1}(\overline{\Omega})$$\cap C^{2,\beta}(\Omega)$, and 
\begin{equation}\label{mariko}
	\lim_{\l'\downarrow 0}\l' u_{\l'}(x)=d,\quad 
	\lim_{\l'\downarrow 0}(u_{\l'}-u_{\l'})(x_0)=u(x)\qquad 
	\hbox{uniformly on} \quad \overline{\Omega}. 
\end{equation} 
 From the usual stability result (\cite{users}n), it is clear that the pair $(d,u)$ satisfies (\ref{1})-(\ref{2}).\\

	As for the uniqueness of the number $d$, we give the following theorem in which 
we consider (\ref{1})-(\ref{2}) in the framework of viscosity solutions. 

\begin{theorem}{\bf Theorem 2.5.} Assume that $\Omega$ is (\ref{doma}), and that (\ref{unif}), (\ref{lips}), (\ref{gamma}) and (\ref{gbd}) hold. Then, the number $d$ such 
that  (\ref{1})-(\ref{2}) has a viscosity solution $u$ is unique.
\end{theorem}

{\it Proof of Theorem 2.5.} We argue by contradiction. Let $(d_1,u_1)$ and $(d_2,u_2)$ be two pairs 
 satisfying (\ref{1})-(\ref{2}) in the sense of viscosity solutions. We assume $d_1>d_2$. 
 First, we show the following Lemma. 

\begin{lemma}{\bf Lemma 2.6.} Let $v=u_1-u_2$. Then, $v$ satisfies
\begin{equation}\label{mm}
-M^{+}(\H v)+\inf_{\a\in A}\{-\sum_{i=1}^n b_i^{\a} \frac{\p v}{\p x_i}\}\le 0 \qquad 
\hbox{in} \quad \Omega,
\end{equation}
\begin{equation}\label{dd}
<\n v, \g>\le d_2-d_1<0 \qquad \hbox{on}\quad \p\Omega,
\end{equation}
where 
\begin{equation}\label{ita}
M^{+}(X)=\sup_{\l_1 I\le A\le \L_1 I} Tr (AX) \qquad X\in {\bf S}^n.
\end{equation}
\end{lemma}
{\it Proof of Lemma 2.6.} Let $\phi\in C^2(\overline{\Omega})$ be such that 
 $u-\phi$ takes its local strict maxixum at $\overline{x}\in \overline{\Omega}$. 
 From the definition of viscosity solutions, we are to show the following.\\

(i) If $\overline{x}\in \Omega$, 
$$
-M^+(\H \phi(\overline{x}))
+\inf_{\a\in A}\{<-b^{\a}(\overline{x}),\phi(\overline{x})>\}\le 0.
$$  

(ii) If $\overline{x}\in \p\Omega$, 
$$
-M^+(\H \phi(\overline{x}))
+\inf_{\a\in A}\{<-b^{\a}(\overline{x}),\phi(\overline{x})>\}\le 0,
$$  
or
$$
<\phi(\overline{x}),\g(\overline{x})>\le d_2-d_1. 
$$

\underline{Step 1.} We shall show (i)  by the contradiction argument. Thus, assume 
\begin{equation}\label{case1}
-M^+(\H \phi(\overline{x}))
+\inf_{\a\in A}\{<-b^{\a}(\overline{x}),\phi(\overline{x})>\}> 0,
\end{equation}
and we shall look for a contradiction. Define, for $\beta>0$
$$
	\Psi_{\beta}(x,y)=u_1(x)-u_2(y)-\phi(\frac{x+y}{2})-\beta|x-y|^2 \qquad 
	\hbox{in}\quad \Omega\times \Omega,
$$
and let $(x_{\beta},y_{\beta})$ be the maximum point of $\Psi_{\beta}$. It is well known  (see \cite{users}n) that 
$$
	(x_{\beta},y_{\beta})\to(\overline{x},\overline{x}), \quad 
	\beta |x_{\beta}-y_{\beta}|^2\to 0\qquad \hbox{as}\quad \beta\to \infty,
$$ 
and that for any $\e>0$, there exist $X$, $Y$ $\in {\bf S}^n$ such that 
$$
(\frac{1}{2}\n \phi(\frac{x_{\beta}+y_{\beta}}{2})+2\beta (x_{\beta}-y_{\beta}),X)
	\in J^{2,+}_{\Omega}u_1(x_{\beta}),
$$
$$
(-\frac{1}{2}\n \phi(\frac{x_{\beta}+y_{\beta}}{2})+2\beta (x_{\beta}-y_{\beta}),Y)
	\in J^{2,-}_{\Omega}u_2(y_{\beta}),
$$
and 
\begin{equation}\label{XY}
-(\frac{1}{\e}+\|A\|)I\le 
	\left( 
\begin{array}{cc}
X & O\\
O & -Y
\end{array}
\right)
\le A+\e A^2,
\end{equation}
where by denoting $\psi(x,y)=\phi(\frac{x+y}{2})+\beta |x-y|^2$, 
$$
	A=D^2\psi(x_{\beta},y_{\beta})\in {\bf S}^{2n}, \quad 
	\|A\|=\sup \{|<A\xi,\xi>|:|\xi|\le 1\}.
$$
Now, by using the definition of viscosity solution for $u_i$ ($i=1,2$), 
$$
F(x_{\beta},\frac{1}{2}\n \phi(\frac{x_{\beta}+y_{\beta}}{2})+2\beta (x_{\beta}-y_{\beta}),X)\le 0,
$$
$$
F(y_{\beta},-\frac{1}{2}\n \phi(\frac{x_{\beta}+y_{\beta}}{2})+2\beta (x_{\beta}-y_{\beta}),Y)\ge 0,
$$
and by taking the differences of two inequalities, using the form of (\ref{hjb}), 
for any small $\d>0$ there exists a control $\a'\in A$ such that
$$
\{-Tr(A^{\a'}(x_{\beta})X)-<\frac{1}{2}\n \phi(\frac{x_{\beta}+y_{\beta}}{2}),b^{\a'}(x_{\beta})>\}\qquad \qquad \qquad \qquad \qquad
$$
\begin{equation}\label{sa}
-
\{-Tr(A^{\a'}(y_{\beta})Y)-<\frac{1}{2}\n \phi(\frac{x_{\beta}+y_{\beta}}{2}),b^{\a'}(y_{\beta})>\}\le \d.
\end{equation}
By taking $\e=\frac{1}{\beta}$ in (\ref{XY}), and multiplying the rightmost inequality in 
(\ref{XY}) by the symmetric matrix
$$
\left(
\begin{array}{cc}
\sigma^{\a'}(x_{\beta})^t\sigma^{\a'}(x_{\beta})&\sigma^{\a'}(y_{\beta})^t\sigma^{\a'}(x_{\beta})\\
\sigma^{\a'}(x_{\beta})^t\sigma^{\a'}(y_{\beta})&\sigma^{\a'}(y_{\beta})^t\sigma^{\a'}(y_{\beta})
\end{array}
\right),
$$
and taking traces, we have 
$$
	Tr (A^{\a'}(x_{\beta})X)-Tr (A^{\a'}(y_{\beta})Y)-Tr (\H \phi(\overline{x}) A^{\a'}(\overline{x}))\le L^2 \beta |x_{\beta}-y_{\beta}|^2+o(\beta^{-1})
$$ 
as $\beta\to \infty$, 
where $L>0$ is the Lipschitz constant in (\ref{lips}) ( or $K$ in (\ref{bdd})). (See \cite{users}, H. Ishii and P.-L. Lions \cite{ipl}n for this techniques.) Therefore from (\ref{sa}), for any $\e>0$ there exists $\a'\in A$ such that 
$$
-Tr(\H \phi(\overline{x})A^{\a}(\overline{x}))-<\n \phi(\overline{x}),b^{\a'}(\overline{x})>\le \d+o(\beta^{-1}), 
$$
which contradicts to (\ref{case1}), since $\d>0$ is arbitrary. Thus, we showed (i).\\

\underline{Step 2.} We shall prove (ii). First of all, from the usual technique to treat 
the Neumann boundary condition in the theory of viscosity solutions, we may replace the 
conditions to 
\begin{equation}\label{for1}
d_1+<\n u_1,\g >-g(x)\le -\d \qquad \hbox{on}\quad \p\Omega,
\end{equation}
\begin{equation}\label{for2}
d_2+<\n u_2,\g >-g(x)\ge \d \qquad \hbox{on}\quad \p\Omega,
\end{equation}
where $\d>0$ is a small number. (See \cite{users}n.) 
 Then, we assume that (ii) does not hold, and shall look for a contradiction. So, let 
\begin{equation}\label{case2}
-M^{+}(\H \phi(\overline{x}))+\inf_{\a\in A}\{
<-b^{\a}(\overline{x}),\n \phi(\overline{x})>\}>0,
\end{equation}
\begin{equation}\label{case2'}
<\n \phi(\overline{x}),\g(\overline{x})>>d_2-d_1.
\end{equation}
 It is well known (\cite{new}n) that since $\p\Omega$ is $C^{3,1}$, by putting 
$$
L(x,y)=\inf\{
\int_0^1 c_{ij}(\xi(t))\dot{\xi}_i\dot{\xi}_j dt\quad |\quad \xi\in C^1([0,1];{\bf R}^n), 
\quad \xi(0)=y,\quad \xi(1)=x 
\},
$$
where $c_{ij}(x)$ is a smooth function, say in $C^3(\overline{\Omega})$ such that for 
${\bf n}=(n_i)_i$
$$
\sum_j c_{ij}(x)\g_j(x)=n_i(x) \qquad \hbox{any}\quad 1\le i\le n,\quad x\in \p\Omega, 
$$
 we have:
\begin{equation}\label{extin}
<\g(x),\n_x L(x,y)><\frac{1}{C}|y-x|^2 \qquad \hbox{any}\quad x\in \p\Omega, \quad y\in \Omega,
\end{equation} 
where $C>0$ is a constant. 
Define, for $\beta>0$
$$
	\Psi_{\beta}(x,y)
=u_1(x)-u_2(y)-\phi(\frac{x+y}{2})-\beta L(x,y)\qquad \qquad \qquad \qquad \qquad
\qquad \qquad \qquad 
$$
$$
	+(d_1-g)<\g(\overline{x}),x-y>+|x-\overline{x}|^4+\frac{1}{2}<\n \phi(\overline{x}),x-y> \qquad 
	\hbox{in}\quad \Omega\times \Omega. 
$$
 Set 
\begin{eqnarray}
&&\psi(x,y)=\phi(\frac{x+y}{2})+\beta L(x,y)
	-(d_1-g)<\g(\overline{x}),x-y>-|x-\overline{x}|^4\nonumber\\
	&&-\frac{1}{2}<\n \phi(\overline{x}),x-y>.\nonumber
\end{eqnarray}
Let $(x_{\beta},y_{\beta})$ be the maximum point of $\Psi_{\beta}$. As in Step1, it is  known  (see \cite{users}) that 
$$
	(x_{\beta},y_{\beta})\to(\overline{x},\overline{x}), \quad 
	\beta |x_{\beta}-y_{\beta}|^2\to 0\qquad \hbox{as}\quad \beta\to \infty,
$$ 
and that for any $\e>0$, there exist $X$, $Y$ $\in {\bf S}^n$ such that 
$$
(\n_x \psi(x_{\beta},y_{\beta}),X)
	\in J^{2,+}_{\Omega}u_1(x_{\beta}),
\quad
(-\n_y \psi(x_{\beta},y_{\beta}),Y)
	\in J^{2,-}_{\Omega}u_2(y_{\beta}),
$$
which satisfy (\ref{XY}) with $A=D^2 \psi\in {\bf S}^{2n}$.  \\
	 If $(x_{\beta},y_{\beta})\in \p\Omega$, by using (\ref{extin})  we calculate
\begin{eqnarray}
&&<\n \psi(x_{\beta},y_{\beta}),\g(x_{\beta})>+d_1-g(x_{\beta})
=<\frac{1}{2}\n \phi(\frac{x_{\beta}+y_{\beta}}{2}),\g(x_{\beta})>\nonumber\\
&&+2\beta<\g(x_{\beta}),\n_x L(x_{\beta},y_{\beta})>-(d_1-g)<\g(x_{\beta}),\g(\overline{x})>
\nonumber\\
&&
-4|x_{\beta}-\overline{x}|^2<\g(x_{\beta}),x_{\beta}-\overline{x}>
-<\g(x_{\beta}),\frac{1}{2}\n \phi(\overline{x})>+d_1-g
\nonumber\\
&&
\ge -\frac{\beta}{C}|x_{\beta}-y_{\beta}|^2+O(|x_{\beta}-z|^3)\ge o(1) 
\qquad \hbox{as} \quad \beta\to \infty.\nonumber
\end{eqnarray}
\begin{eqnarray}
&&<-\n \psi(x_{\beta},y_{\beta}),\g(y_{\beta})>+d_2-g(y_{\beta})
=<-\frac{1}{2}\n \phi(\frac{x_{\beta}+y_{\beta}}{2}),\g(y_{\beta})>\nonumber\\
&&-2\beta<\g(y_{\beta}),\n_y L(x_{\beta},y_{\beta})>-(d_1-g)<\g(y_{\beta}),\g(\overline{x})>
\nonumber\\
&&
+<\g(y_{\beta}),\frac{1}{2}\n \phi(\overline{x})>+d_2-g
\nonumber\\
&&
\le \frac{\beta}{C}|x_{\beta}-y_{\beta}|^2+d_2-d_1+o(1)\le o(1) 
\qquad \hbox{as} \quad \beta\to \infty.\nonumber
\end{eqnarray}
(In the last inequality, we used the assumption $d_1>d_2$.)\\
	Therefore, by taking account of (\ref{for1}) and (\ref{for2}), regardless the fact that $x_{\beta}$, $y_{\beta}$ $\in \Omega$ or $\in \p\Omega$, we have the following.
$$
F(x_{\beta},\n \psi(x_{\beta},y_{\beta}),X)\le o(1)\qquad \hbox{as}\quad \beta\to \infty,
$$
$$
F(y_{\beta},-\n \psi(x_{\beta},y_{\beta}),Y)\ge o(1)\qquad \hbox{as}\quad \beta\to \infty.
$$
 The rest of the argument to obtain a contradiction from the above two inequalities 
 is similar to that of Step 1, and we omit it here.
 
	Now, we go back to the proof of Theorem 2.5, which is immediate from Lemma 2.6. 
 From the strong maximum principle (Lemma A), $v$, which is not constant, attains its maximum at some point $x_1\in \p\Omega$
$$
	v(x_1)>v(x) \qquad \hbox{any}\quad x\in \Omega.
$$
However, as we have seen in the proof of Theorem 2.1 in Step 1, this is not 
compatible with $<\n v,\g>\le d_2-d_1$ on $\p\Omega$, in the sense of viscosity solutions. Thus, we have proved $d_1=d_2$ must be hold. \\

	If we consider the uniqueness of $d$ in the framework of the $C^{1,1}(\overline{\Omega})$ solutions, the proof is much simpler. We add this as follows.
\\

\begin{theorem}{\bf Proposition 2.7.} Assume that $\Omega$ is (\ref{doma}), and that 
(\ref{unif}), (\ref{lips}) and (\ref{gamma}) hold. Moreover, assume that $F$ satisfies the 
following comparison: for a subsolution $u$ and a supersolution $v$ of (\ref{1}) such that 
$u\le v$ on $\p\Omega$,  $u\le v$ in $\overline{\Omega}$.
Then,  
 the number $d$ such that (\ref{1})-(\ref{2}) has a solution $u\in C^{1,1}(\overline{\Omega})$ is unique.
\end{theorem}

{\it Proof of Proposition 2.7.} We assume that there are two pairs $(d_1,u_1)$ and 
$(d_2,u_2)$ which satisfy (\ref{1})-(\ref{2}) such that $d_1>d_2$ and 
$u_i\in C^{1,1}(\overline{\Omega})$ ($i=1,2$). By adding a constant if necessary, 
we may assume that there is a point $x_0\in \p\Omega$ such that $u_1(x_0)=u_2(x_0)$ and 
$$
	u_1(x)\le u_2(x) \qquad \hbox{on}\quad \p\Omega.
$$
Put $v=u_2-u_1$, which satisfies 
$$
	<\n v(x),\g(x)>=d_1-d_2>0,\quad v(x)\ge 0 \qquad \hbox{on}\quad \p\Omega.
$$
From the comparison for (\ref{1}), 
$$
v(x)\ge 0 \qquad \hbox{any}\quad x\in \overline{\Omega}.
$$
However, at $x_0\in \p\Omega$, $v(x_0)=0$ and $<\n v(x_0),\g(x_0)>>0$ in the classical sense. 
Thus, we get a contradiction and $d_1=d_2$. \\

\section{Long time averaged reflection force in  half spaces.}

	In this section, the existence and uniqueness of the number $d$ in 
(\ref{1})-(\ref{2}) is shown in the case that $\Omega$ satisfies (\ref{dom2}), with  a supplement 
boundary condition at $x_n=\infty$. We denote
$$
\Omega=\{(x',x_n)|\quad x_n\ge f(x'),\quad x'\in ({\bf R}\backslash{\bf Z})^{n-1}\},
$$
$$
\G_0=\p\Omega=\{(x',x_n)|\quad x_n=f(x'), \quad x'\in ({\bf R}\backslash {\bf Z})^{n-1}\},
$$

\newpage
where $f(x')$ is periodic in $x'\in(\R\backslash {\bf Z})^{n-1}$ and is $C^{3,1}$. 
 Our goal is  to find a unique number $d$ which admits a viscosity solution $u$ of (\ref{1})-(\ref{2}) such that 
\begin{equation}\label{mugen}
u \quad \hbox{is bounded and periodic in $x'$.}
\end{equation}
We begin with the uniqueness of $d$.

\begin{theorem}{\bf Theorem 3.1.} Assume that $\Omega$ is (\ref{dom2}), and that (\ref{unif}), (\ref{lips}), (\ref{gamma}) and (\ref{gbd}) hold. 
 Moreover, assume that 
\begin{equation}\label{drift}
b^{\a}_n(x)\le 0 \qquad \hbox{any}\quad x\in \Omega,\quad \a\in A.
\end{equation}
Then, the number $d$ 
such that (\ref{1})-(\ref{2}) and (\ref{mugen}) has a viscosity solution $u$ is unique.
\end{theorem}

{\it Proof of Theorem 3.1.} We argue by contradiction. Assume that there exist two pairs 
 $(d_1,u_1)$ and $(d_2,u_2)$ which satisfy (\ref{1})-(\ref{2}) and (\ref{mugen}), and that 
 $d_1>d_2$. By using a similar argument to the proof of Lemma 2.6, $v=u_1-u_2$ 
is a subsolution of 
\begin{equation}\label{pucci}
-M^{+}(\H v)+\inf_{\a}\{<-b^{\a}(x),\n v>\}\le 0 \qquad \hbox{in}\quad \Omega,
\end{equation}
\begin{equation}\label{pucci2}
<\n v,\g(x)>=d_2-d_1<0 \qquad \hbox{on} \quad \p\Omega,
\end{equation}
where $M^{+}$ is the Pucci operator defined in (\ref{ita}) (See \cite{swi}n).
For $R>0$ large enough, let 
$$
\Omega_R=\{(x',x_n)|\quad f(x')\le x_n\le R\},
$$
and define
$$
M_R=\sup_{\overline{\Omega_R}}|v|.
$$ 
(Remark that $v$ is periodic in $x'\in ({\bf R}\backslash {\bf Z})^{n-1}$ and the above 
 supremum is well-definded.)
Let $x_0\in \G_0$ be a point such that $v(x_0)=\sup_{x\in\G_0}v(x)\equiv M_0$. 
  Let $(x_c',c)\in \G_0$ be a point such that 
$$
c\le x_n\qquad \hbox{any} \quad (x',x_n)\in \G_0. 
$$
We take 
\begin{equation}\label{wR}
w_{R}(x',x_n)\equiv \frac{M_R-M_0}{R-c}(x_n-c)+M_0 \qquad (x',x_n)\in \Omega.
\end{equation} 
 Since $\frac{M_R-M_0}{R-c}\ge 0$, from (\ref{drift}) 
$$
-M^{+}(\H w_R)+\inf_{\a}\{<-b^{\a}(x),\n w_R>\}\ge 0 \qquad \hbox{in}\quad \Omega_R,
$$
$$
w_{R|\G_0}=\frac{M_R-M_0}{R-c}(x_n-c)+M_0\ge M_0, 
$$
$$
w_{R|\G_R}=M_R.
$$
 Thus, by using the comparison argument, we get
$$
v\le w_R \qquad \hbox{in}\quad \overline{\Omega_R},\quad \hbox{any}\quad R>0\quad 
\hbox{large enough.}
$$
 By (\ref{mugen}), tending $R\to \infty$, this yields
$$
v\le M_0 \qquad \hbox{in}\quad \Omega.
$$
Therefore, $v$ takes its maximum on $\G_0$. Finally, by using the strong maximum principle 
 (Lemma A), (\ref{pucci}) and (\ref{pucci2}) yields a contradiction as we argued in 
 the proof of Theorem 2.1,  Step1. Thus, $d_1=d_2$ must hold.\\

{\bf Remark 3.1.} (Counter example.) If we do not assume the boundary condition at 
infinity (\ref{mugen}), $d$ is not unique in general. For example, consider 
\begin{equation}\label{exe}
-\Delta u=0 \qquad \hbox{in} \quad \{x_n\ge 0\}\subset \Rn,
\end{equation}
\begin{equation}\label{exb}
d+<\n u,{\bf n}(x)>=0 \qquad \hbox{on} \quad \{x_n=0\}\subset \Rn,
\end{equation}
where ${\bf n}$ is the outward unit normal, and the solution $u$ is periodic in 
$x'$$=(x_1,...,x_{n-1})$. Then, for any $c$, $d\in R$, $u=-dx_n+c$ is the solution of 
 (\ref{exe})-(\ref{exb}). Thus, the number $d$ in (\ref{exb}) is not unique. 
 (Green's first identity does not hold in the half space.) \\

	Next, for the existence of $d$ we approximate (\ref{1})-(\ref{2}) and (\ref{mugen}) by
$$
F(x,\n u_{\l}^R,\H u_{\l}^R)=0 \qquad \hbox{in} \quad \Omega_R =\{(x',x_n)|\quad f(x')\le x_n\le R\},
$$ 
\begin{equation}\label{pR}
<\n u_{\l}^R,{\bf n}(x)>=0 \qquad \hbox{on} \quad \G_R=\{(x',x_n)|\quad x_n=R\},
\end{equation}
$$
\l u_{\l}^R+<\n u_{\l}^R,\g(x)>-g(x)=0 \qquad \hbox{on}\quad \p\Omega=\G_0=
\{x_n=f(x')\},
$$
where $R>0$ is large enough so that $\G_R$ and $\G_0$ do not intersect, say $R\ge R_0$.
 We examine the regularity of $u_{\l}^R$ uniformly in $\l\in(0,1)$ and $R>R_0$. 
 
\begin{proposition}{\bf Proposition 3.2.} Assume that $\Omega$ is (\ref{dom2}), and that 
(\ref{unif}), (\ref{gamma}), (\ref{bdd}) and (\ref{gbd}) hold.  Let $R>R_0$ be fixed, and let $u_{\l}^R$ be the solution of (\ref{pR}). Then, there exists a number $d_R$ and a function $u_{R}$ such that 
\begin{eqnarray}\label{Rlim}
\lim_{\l\downarrow 0} \l u_{\l}^R(x)&=& d_R,\nonumber \\
\lim_{\l'\downarrow 0} (u_{\l'}^R(x)-u_{\l'}^R(x_0))&=& u_R(x) \qquad \hbox{uniformly in}\quad \overline{\Omega_R},
\end{eqnarray}
where $\l'\to 0$ is a subsequence of $\l\to 0$, and $x_0$ is an arbitrarily fixed point in $\Omega_{R_0}$. The pair $(d_R,u_R)$ satisfies 
$$
F(x,\n u_R,\H u_R)=0 \qquad \hbox{in} \quad \Omega_R,
$$ 
\begin{equation}\label{R}
<\n u_R,{\bf n}(x)>=0 \qquad \hbox{on} \quad \G_R,
\end{equation}
$$
d_R+<\n u_R,\g(x)>-g(x)=0 \qquad \hbox{on}\quad \p\Omega=\G_0.
$$
The number $d_R$ is the unique number such that (\ref{R}) has a viscosity solution. 
Moreover, there exists a constant $M>0$ such that 
\begin{equation}\label{Rb}
|u_R-u_R(x_0)|_{L^{\infty}(\overline{\Omega_R})}<M \qquad \hbox{any}\quad  R>R_0,
\end{equation}
\begin{equation}\label{Rbd}
|\n u_R|_{L^{\infty}(\overline{\Omega_{R}})}<M \qquad \hbox{any}\quad R>R_0.
\end{equation}
\end{proposition}
{\it Proof of Proposition 3.2.}  We devide the proof into three steps.\\ 
\underline{Step 1.} First, we shall see
\begin{equation}\label{uRb}
|u_{\l}^R(x)-u_{\l}^R(x_0)|\le M \qquad \hbox{any}\quad \l\in (0,1),\quad  R>R_0.
\end{equation}
So, put $v_R=u_R-u_R(x_0)$. 
Assume that 
$$
(\e_{\l}^R)^{-1}\equiv |v_{\l}^R|_{L^{\infty}(\overline{\Omega_R})}\to \infty \qquad \hbox{as} \quad \l\to 0,
\quad R\to \infty,
$$
and we seek a contradiction. Put $w_{\l}^{R}\equiv \e_{\l}^Rv_{\l}^R$ which satisfies 
$$
F(x,\n w_{\l}^R,\H w_{\l}^R)=0 \qquad \hbox{in} \quad \Omega_R,
$$
$$
<\n w_{\l}^R,{\bf n}(x)>=0 \qquad \hbox{on} \quad \G_R,
$$
$$
<\n w_{\l}^R,\g(x)>=\e_{\l}^R (g-\l u_{\l}^R) \qquad \hbox{on} \quad \G_0.
$$
Since $|w_{\l}^R|_{L^{\infty}(\overline{\Omega_R})}=1$ ($w_{\l}^R(x_0)=0$), 
$$
w^{\ast}(x)=\limsup_{R\to \infty, \l\downarrow 0, y\to x} w_{\l}^R(y),\quad w_{\ast}(x)=\liminf_{R\to \infty, \l\downarrow 0, y\to x} w_{\l}^R(y),
$$
are well-definded. From the uniform ellipticity (\ref{unif}) and the Krylov-Safonov interior estimate,  for any $V\subset\subset \Omega$ there exists a constant 
 $M_V>0$ such that 
$$
|\n w_{\l}^R|_{L^{\infty}(V)}\le M_V \qquad \hbox{any}\quad \l\in(0,1), \quad 
R>R_0.
$$
Thus, since $w_{\lambda}(x_0)=0$ $(\forall \lambda\in (0,1))$,
\begin{equation}\label{ks}
w^{\ast}(x_0)=w_{\ast}(x_0)=0.
\end{equation}
Moreover from the strong maximum principle (Lemma A), for any $R>R_0$ and $\l\in (0,1)$, 
$w_{\l}^{R}$ must take its maximum and minimum on $\G_0$. (If it takes a maximun or a minimum on $\G_R$, we have a 
contradiction to $<\n w_{\l}^R,{\bf n}(x)>=0$ ($x\in \G_R$) in the sense of viscosity solutions as we have seen in the proof of Theorem 2.1, Step 1.) Hence,
\begin{equation}\label{kss}
 |w^{\ast}|_{L^{\infty}(\overline{\Omega_R})}=1 \quad \hbox{or}\quad |w_{\ast}|_{L^{\infty}(\overline{\Omega_R})}=1\qquad \hbox{any}\quad R>R_0.
\end{equation}
 Hereafter, we assume that $|w^{\ast}|_{L^{\infty}(\overline{\Omega_R})}=1$. (The case of $|w_{\ast}|_{L^{\infty}(\overline{\Omega_R})}=1$ can be treated similarly.)
The upper semicontinuous function $w^{\ast}$ is a viscosity solution of 
\begin{equation}\label{ws}
F(x,\n w^{\ast},\H w^{\ast})\le 0 \qquad \hbox{in} \quad \Omega,
\end{equation}
\begin{equation}\label{wss}
<\n w^{\ast},\g(x)>\le 0 \qquad \hbox{on} \quad \G_0.
\end{equation}

We remark that $w^{\ast}$ takes its maximum on $\G_0$, as $w_{\lambda}^R$ ($R>R_0$, $\lambda\in (0,1)$) does so. ($w^{\ast}$ is periodic in $x'\in ({\bf R}\backslash {\bf Z})^{n-1}$.) Then, 
by the strong maximum principle (Lemma A) and the fact that $w^{\ast}$ is not constant ((\ref{ks}), (\ref{kss})), 
 (\ref{ws})-(\ref{wss}) lead a contradiction. (See the proof of Theorem 2.1, Step 1.) Therefore, there exists a constant $M>0$ such 
that 
$$
|u_{\l}^R(x)-u_{\l}^R(x_0)|\le M \qquad \hbox{any}\quad \l\in (0,1),\quad  R>R_0.
$$

\underline{Step 2.} Next, we shall show (\ref{Rlim}) and (\ref{Rbd}).
 For  this purpose, we are to have the a priori estimates of $|\n u_{\l}^R|$ and $|\H u_{\l}^R|$.  Put
\begin{equation}\label{wsss}
w_{\l}^R=\frac{u_{\l}^R-u_{\l}^R(x_0)}{|\n (u_{\l}^R-u_{\l}^R(x_0))|_{L^{\infty}(\overline{\Omega_{R}})}}. 
\end{equation}
Remark that $w_{\l}^R$ is a solution of 
$$
F(x,\n w_{\l}^R,\H w_{\l}^R)=0 \qquad \hbox{in} \quad \Omega_R,
$$
\begin{equation}\label{nR}
<\n w_{\l}^R,{\bf n}(x)>=0 \qquad \hbox{on} \quad \G_R,
\end{equation}
\begin{equation}\label{n0}
\l w_{\l}^R+<\n w_{\l}^R,\g(x)>-\overline{g}=0 \qquad \hbox{on} \quad \G_0,
\end{equation}
where 
$$
\overline{g}=\frac{g}{|\n (u_{\l}^R-u_{\l}^R(x_0))|_{L^{\infty}(\overline{\Omega_{R}})}}.
$$
Taking account of the periodicity in $x_i$ $(i=1,...,n-1)$, the above problem is reduced 
to the case of bounded domains treated in $\S$ 2. Despite the existence of the different boundary condition (\ref{nR}) on $\G_R$, the argument in $\S$ 2 (and \cite{lt}n) works with a minor modification. (We do not rewrite it here.) Thus, the a priori estimate: 
$$
|\H w_{\l}^R|_{L^{\infty}(\overline{\Omega_R})}\le M \qquad \hbox{any}\quad \l\in (0,1),\quad  R>R_0,
$$
where $M>0$ is a constant, which leads
\begin{equation}\label{est}
|\H u_{\l}^R|_{L^{\infty}(\Omega_R)}\le M(|\n u_{\l}^R|_{L^{\infty}(\Omega_R)}+1) 
\qquad \hbox{any}\quad \l\in (0,1),\quad  R>R_0.
\end{equation}
As in $\S$ 2, we use the interpolation inequality in Lemma B, 
 with the function $u_{\l}^R-u_{\l}^{R}(x_0)$, $D=\Omega_{R}$, $j=1$, $k=2$ and  $\alpha =\beta=0$.  That is, the interpolation inequality becomes: 
\begin{equation}\label{fine}
|\n u_{\l}^R|_{L^{\infty}(\overline{\Omega_{R}})}\le C_{\e}|u_{\l}^R-u_{\l}^R(x_0)|_{L^{\infty}(\overline{\Omega_{R}})}+\e|\H u_{\l}^R|_{L^{\infty}(\overline{\Omega_{R}})}.
\end{equation}
 By combining (\ref{Rb}), (\ref{est}) and (\ref{fine}),  
$$
|\H u_{\l}^R|_{L^{\infty}(\Omega_{R})}\le M \qquad \hbox{any} \quad \l\in (0,1),\quad R>R_0,
$$
$$
|\n u_{\l}^R|_{L^{\infty}(\Omega_{R})}\le M \qquad \hbox{any} \quad \l\in (0,1),\quad R>R_0.
$$
Thus, by extracting a subsequence $\l'\downarrow 0$, there exists a number 
$d_R$ and a function $u_R$ such that 
$$
\l' u_{\l'}^R \to d_R, \qquad u_{\l'}^R-u_{\l'}^R(x_0)\to u_R, 
$$
and 
$$
|\n u_R|_{L^{\infty}(\overline{\Omega_{R}})}\le M \qquad \hbox{any} \quad R>R_0.
$$
Thus, we  proved (\ref{Rlim}) and (\ref{Rbd}).\\

\underline{Step 3.}
We shall complete the proof  by showing that the above limit $d_R$ is the unique number such that (\ref{R}) has a viscosity solution (and is independent of the choice of $\l'\to 0$).   We argue by contradiction, and 
assume that there exist two pairs $(d_R,u_R)$ and $(d_R',u_R')$ $(d_R>d_R')$ satisfying (\ref{R}). Denote $v=u_R-u_R'$.  A similar argument used in the proof of Lemma 2.6 leads  
$$
-M^{+}(\H v)+\inf_{\a\in A}\{<-b^{\a}(x),\n v>\}\le 0 \qquad \hbox{in}\quad \Omega_R,
$$
$$
<\n v,{\bf n}(x)>\le 0 \qquad \hbox{on}\quad \G_R,
$$
$$
<\n v,\g(x)>\le d_R'-d_R\qquad \hbox{on}\quad \G_0. 
$$
Since $v$ is not constant, from the strong maximum principle (Lemma A), $v$ attains its 
 maximum at $x_0\in\G_0$: 
$$
v(x_0)>v(x) \qquad \hbox{any}\quad x\in \Omega_R.
$$
 However, as we have seen in the proof of Theorem 2.1 Step1, since $d_R'-d_R<0$, it is not compatible with the preceding boundary conditions on $\G_0$ and $\G_R$. Therefore, we get a contradiction and $d_R=d'_R$ must hold. \\

\begin{theorem}{\bf Theorem 3.3.} Assume that $\Omega$ is (\ref{dom2}), and that (\ref{unif}), (\ref{gamma}), (\ref{bdd}) and (\ref{gbd}) hold. Then, 
there exists a unique number $d$ such that  (\ref{1})-(\ref{2}) and (\ref{mugen}) has a viscosity solution $u$.
\end{theorem}
{\it Proof of Theorem 3.3.} By comparison, there exists a constant $C>0$ such that 
$$
|\l u_{\l}^R|_{L^{\infty}(\overline{\Omega_R})}\le C
\qquad \hbox{any}\quad  \l\in (0,1),\quad  R>R_0, 
$$ 
and  thus $|d_R|<C$ for any $R>R_0$.  Therefore, by using (\ref{Rb}) and (\ref{Rbd}), we can extract a subsequence $R'\to \infty$ such that there exist a number $d$ and a function $u$ such that
$$
d_{R'}\to d  \qquad \hbox{as} \quad R'\to \infty,
$$
$$
u_{R'}\to u \qquad \hbox{as}\quad R'\to \infty, \quad \hbox{locally uniformly in}\quad 
\overline{\Omega}. 
$$
 From the stability results,
$$
F(x,\n u,\H u)=0 \qquad \hbox{in}\quad \Omega,
$$
$$
d+<\n u,\g(x)>-g(x)=0 \qquad \hbox{on} \quad \G_0,
$$
$$
|u|_{L^{\infty}(\overline{\Omega})}<M.
$$
The uniqueness of $d$ was proved in Theorem 3.1, and we can end the proof.\\

{\bf Remark 3.2.} From the view point of the stochastic process (\ref{process}), the 
approximating system (\ref{R}) gives a kind of boundary condition at infinity. It forces 
 the admissible trajectories of (\ref{process}) (corresponding to (\ref{1})-(\ref{2}) and (\ref{mugen})) to be pushed back inward at some finite $x_n=R$. Therefore, the condition (\ref{drift}) is quite reasonable. (In \cite{ben}n, the ergodic problem in unbounded 
 domain (not on the boundary like (\ref{2})) is solved with the condition $\lim_{x\to \infty} b_n^{\a}(x)=-\infty$, which 
is stronger than (\ref{drift}).) \\
\section{Remarks on some degenerate cases.}

	The number $d$ in (\ref{1})-(\ref{2}) exists even for degenerate operators. 
 In this section, we 
give a sufficient condition for the existence (in a weeker sense) and two  classes of operators  satisfying the sufficient condition. The following two examples illustrate 
the existence and non-uniqueness of $d$. In the case of degenerate operators, the uniqueness does not hold in general. \\ 

\begin{example}{\bf Example 4.1.} Consider
$$
|\n u|=0\qquad \hbox{in}\quad \Omega,
$$
\begin{equation}\label{dexp}
d+<\n u,{\bf n}(x)>-g(x)=0 \qquad \hbox{on}\quad \p\Omega,
\end{equation}
where $\Omega\subset \Rn$  is a bounded open domain with a smooth boundary $\p\Omega$, ${\bf n}$ is the outward unit normal to $\Omega$, and $g$ is Lipschitz continuous on $\p\Omega$.  
Then,  any $d$ such that
$$
d\le \min_{x\in \p\Omega} g(x)
$$
 and $u\equiv\hbox{C}$ (constant) 
satisfies (\ref{dexp}) in the sense of viscosity solutions. In fact, 
it is clear that $u$ satisfies the equation in $\Omega$. To see the boundary condition 
in the viscosity sense, 
$$
\max\{|\n u|,d+<\n u,{\bf n}(x)>-g(x)\}\ge 0 \qquad \hbox{on} \quad \p\Omega,
$$
shows that $u$ is a supersolution on $\p\Omega$. For any $\phi\in C^1$ such that $u-\phi$ 
takes its strict maximum at $x_0\in \p\Omega$,  if  $d\le \min_{\p\Omega} g$ then 
$$
<\n\phi,{\bf n}(x)>\quad \le 0\quad \le g(x)-d \qquad \hbox{on}\quad  \p\Omega. 
$$
Thus,
$$
\min \{|\n u|,d+<\n u,{\bf n}(x)>-g(x)\}\le 0 \qquad \hbox{on}\quad  \p\Omega,
$$
 in the sense of viscosity solutions, and $u$ is a subsolution on $\p\Omega$.
\end{example}
\begin{example}{\bf Example 4.2.}
Let $\Omega$$=(\R/{\bf Z})\times (0,1)$$\subset \R^2$ (periodic in $x_1$). Consider 
$$
-\frac{\p^2 u}{\p x_1^2}+|\frac{\p u}{\p x_2}|=0 \qquad \hbox{in}\quad \Omega,
$$
\begin{equation}\label{dexp2}\quad
\end{equation}
$$
d+<\n u, {\bf n}(x)>-g(x)=0 \qquad \hbox{on}\quad \p\Omega,
$$
where ${\bf n}$ is the outward unit normal to $\Omega$, $g$ is Lipschitz continuous on $\p\Omega$. Then, any $d$ such that
$$
d\le \min_{x\in \p\Omega} g(x)
$$
and $u\equiv \hbox{C}$ (constant) 
satisfies (\ref{dexp2}) in the sense of viscosity solutions.
 In fact clearly, $u$ is a viscosity solution in $\Omega$. To see that $u$ is a supersolution on 
$\p\Omega$, suppose for $\phi\in C^1$, $u-\phi$ takes its strict minimum at $x_0\in \p\Omega$. Since $u=C$ on $x_1=0,1$,  
 we remark that such $\phi\in C^2$ must not satisfy $-\frac{\p^2 \phi}{\p x_1^2}(x_0)\le 0$. 
Thus, 
$$
-\frac{\p^2 \phi}{\p x_1^2}(x_0)+|\frac{\p \phi}{\p x_2}(x_0)|\ge 0,
$$
 and $u$ is a viscosity  super solution on $\p\Omega$. The fact that $u$ is a 
 subsolution on $\p\Omega$ is same to  Example 4.1.
\end{example}

{\bf Remark 4.1.} In the above examples the numbers $d$ are not unique. \\

	The operators $F$ studied here are given in (\ref{hjb}) with 
 degenerate coefficients. For such operators, we
 approximate (\ref{1})-(\ref{2}) by
\begin{equation}\label{ep}
-\e\Delta u_{\e}+F(x,\n u_{\e}, \H u_{\e})=0 \qquad \hbox{in}\quad \Omega,
\end{equation}
\begin{equation}\label{eb}
d_{\e}+<\n u_{\e},\g(x)>-g(x)=0 \qquad \hbox{on}\quad \p\Omega,
\end{equation}
where $\e\in (0,1)$. 
 The domain $\Omega$ is either (\ref{doma}) or (\ref{dom2}), and in the case of (\ref{dom2}) the condition at infinity (\ref{mugen}) is added. 
 For any $\e>0$, the existence and the uniqueness of $d_{\e}$ and the existence of $u_{\e}$ come from  Theorems 2.4, 2.5, and 3.3, for (\ref{ep}) is uniformly elliptic.

\begin{proposition}{\bf Proposition 4.1.} Let $\Omega$ be a domain either (\ref{doma}) or (\ref{dom2}). In the case of (\ref{doma}), assume all conditions but (\ref{unif}) in Theorems 2.4 and 2.5.  In the case of (\ref{dom2}), assume all conditions but (\ref{unif}) in Theorem 3.3. (Thus, $F$ is possibly degenerate.)  Let $d_{\e}$ ($\e>0$) be the number such that (\ref{ep})-(\ref{eb}) 
 (and (\ref{mugen}) in the case of (\ref{dom2})) has a viscosity solution $u_{\e}$. Assume that there is a number $M>0$ such that 
\begin{equation}\label{suff}
|u_{\e}-u_{\e}(x_0)|_{L^{\infty}(\Omega)}<M \qquad \hbox{any}\quad \e\in(0,1). 
\end{equation}
Then, there exists a number $d$ (not necessarily unique) such that (\ref{1})-(\ref{2})
 (and (\ref{mugen}) in the case of (\ref{dom2})) has a viscosity subsolution $\underline{u}$ and 
a supersolution $\overline{u}$. 
\end{proposition} 

{\it Proof of Proposition 4.1.} Put $v_{\e}=u_{\e}-u_{\e}(x_0)$. Since $d_{\e}$ is bounded in $\e\in (0,1)$, we can take a subsequence $\e'\to 0$ such that  $\lim_{\e'\to 0}d_{\e}= d$ holds for a constant $d$. From (\ref{suff}),  
$$
v^{\ast}(x)=\limsup_{\e '\downarrow 0,y\to x} v_{\e}(y),\quad 
v_{\ast}(x)=\liminf_{\e '\downarrow 0,y\to x} v_{\e}(y) 
$$
are well-definded. 
 Then, from the usual stability result (see \cite{users}n),  $(d,v^{\ast})$ and  $(d,v_{\ast})$
are respectively viscosity sub and super solutions of (\ref{1})-(\ref{2})
 (and (\ref{mugen}) in case of (\ref{dom2})).\\

\bigskip

{\bf Remark 4.2.} In the above proposition $v^{\ast}\neq v_{\ast}$ in general, and thus 
 the result is weaker than uniformly elliptic cases. \\

\bigskip

	Next, we give a class of operators satisfying (\ref{suff}).
 The first class admits the existence of the uniformly elliptic part: \\

there exists a point $x_0\in \Omega$ such that in a small neighborhood $B(x_0,r)\subset \Omega$ ($r>0$), there exist constants $\l_2$ and $\L_2$ such that $0<\l_2\le \L_2$ and 
\begin{equation}\label{hole}
	\l_2 I\le (a_{ij}^{\a})_{1\le i,j\le n}\le \L_2
	\qquad \hbox{any}\quad \a\in A,\quad x\in B(x_0,r).
\end{equation}

 The second class admits the existence of the "controllability" part (see \cite{ar2}):\\

there exists a point $x_0\in \Omega$ such that for a small neighborhood $B(x_0,r)\subset 
\Omega$ ($r>0$), 
\begin{equation}\label{hole2}
\lim_{|p|\to \infty} F(x,p,X)\to \infty \qquad \hbox{uniformly in}\quad 
x\in \Omega,\quad X\in {\bf S}^n.
\end{equation}

\begin{theorem}{\bf Thorem 4.2.} Let $\Omega$ be a domain either (\ref{doma}) or (\ref{dom2}).  In the case of (\ref{doma}), assume all conditions but (\ref{unif}) in Theorems 2.4 and 2.5.  In the case of (\ref{dom2}), assume all conditions but (\ref{unif}) in Theorem 3.3. (Thus, $F$ is possibly degenerate.) 
Assume also that $F$ satisfies (\ref{bardi1}), (\ref{bardi2}) and (\ref{bardi3}), and that 
 either (\ref{hole}) or (\ref{hole2}) holds. 
 Then, the solutions $u_{\e}$ ($\e>0$) of (\ref{ep})-(\ref{eb}) (and (\ref{mugen}) in the case of  (\ref{dom2}))  satisfy (\ref{suff}).  Moreover, there exists a number $d$ (not necessarily unique) such that (\ref{1})-(\ref{2})
 (and (\ref{mugen}) in the case of (\ref{dom2})) has a viscosity subsolution $\underline{u}$ and 
a supersolution $\overline{u}$. 
\end{theorem}

{\it Proof of Theorem 4.2.} Assume that $(\ref{suff})$ does not hold, and we shall look for  a contradiction. Let $x_0$ be a point satisfying (\ref{hole}) or (\ref{hole2}), and assume that 
$|u_{\e}-u_{\e}(x_0)|_{L^{\infty}(\Omega)}\to \infty$ as $\e>0$ goes to $0$. Put 
$$
v_{\e}=\frac{u_{\e}-u_{\e}(x_0)}{|u_{\e}-u_{\e}(x_0)|_{L^{\infty}(\Omega)}}. 
$$
 The function $v_{\e}$ satisfies 
$$
-\e\Delta v_{\e}+F(x,\n v_{\e},\H v_{\e})=0 \qquad \hbox{in}\quad \Omega,
$$
$$
<\n v_{\e},\g>=\frac{g(x)-d_{\e}}{|u_{\e}-u_{\e}(x_0)|_{L^{\infty}(\Omega)}} \qquad \hbox{on}\quad \p\Omega. 
$$
Since $|v_{\e}|_{L^{\infty}(\Omega)}=1$, 
$$
v^{\ast}(x)=\limsup_{\e\downarrow  0,y\to x}v_{\e}(y),\quad v_{\ast}(x)=\liminf_{\e\downarrow 0, y\to x}v_{\e}(y),
$$
are well definded.  Now, in the case of (\ref{hole}),  we use the Krylov-Safonov inequality
 as before to have
\begin{equation}\label{zero}
v^{\ast}(x_0)=v_{\ast}(x_0)=0. 
\end{equation}
 In  the case of (\ref{hole2}), by using the argument in \cite{ipl},
  \cite{new}n we  have 
 also the uniform continuity of $u_{\e}$ $(\e\in (0,1))$ in $B(r,x_0)$, and   (\ref{zero}) holds. 
 In conclusion,  (\ref{zero}) holds in both cases of (\ref{hole}) and (\ref{hole2}).\\
	We continue the proof, and see easily either 
  $|v^{\ast}|_{L^{\infty}(\overline{\Omega)}}=1$  or $|v_{\ast}|_{L^{\infty}(\overline{\Omega)}}=1$ holds. If $|v^{\ast}|_{L^{\infty}(\overline{\Omega)}}=1$, since 
$$
F(x,\n v^{\ast},\H v^{\ast})\le 0 \qquad \hbox{in}\quad \Omega,
$$
$$
<\n v^{\ast},\g>\le 0\qquad \hbox{on}\quad \p\Omega,
$$
 the strong maximum principle (Lemma A) leads a contradiction, for $v^{\ast}$ is not constant (\ref{zero}). (See the proof of Theorem 2.1, Step 1.)
  If $|v_{\ast}|_{L^{\infty}(\Omega)}=1$, the same argument works, too. Therefore,   $u_{\e}$ satisfies (\ref{suff}), and  Proposition 4.1 leads the remained claim. \\

	As for the uniqueness of $d$, we do not have the general result, and shall 
 give the following Example in which the uniqueness holds.  \\

\begin{example}{\bf Example 4.3.} Let $\Omega$$=\{(x_1,x_2)|\quad x_1\in \R\backslash {\bf Z}, x_2> 0\}$$\subset \R^2$ (periodic in $x_1$). 
Assume that there exists a number $d$ such that 
$$
-\frac{\p^2 u}{\p x_2^2}-\frac{\p u}{\p x_1}=0 \qquad \hbox{in} \quad \Omega,
$$
$$
d+<\n u,{\bf n}(x)>-g(x)=0\qquad \hbox{on} \quad \p\Omega,
$$
where $u$ is bounded, and ${\bf n}$ is the outward unit normal to $\Omega$. 
Then, $d=\int_0^1g(x_1,0)dx_1$.\\

	In fact, by integrating the above problem in $x_1\in [0,1]$,   $\overline{u}(x_2)=\int_0^1 u(x_1,x_2)dx_1$ satisfies  
$$
-\frac{\p^2 \overline{u}(x_2)}{\p x_2^2}=0 \qquad \hbox{in} \quad (0,\infty),
$$
$$
d-\frac{\p \overline{u}(0)}{\p x_2}-\int_0^1 g(x_1,0) dx_1=0 \qquad \hbox{on} \quad x_2=0,
$$
and $\overline{u}$ is bounded. From Theorem 3.3, we know that such a number $d$ is unique. Since $d=\int_0^1 g(x_1) dx_1$ and $\overline{u}\equiv C$ (constant) satisfy the above, we proved the claim.
\end{example}

\section{Homogenization of oscillating Neumann type boundary conditions.}

	In this section, we study the following homogenization problem.
\begin{equation}\label{Go}
G(x,\n u_{\e},\H u_{\e})=\sup_{\a\in \CA}\{-\sum_{ij=1}^2 a_{ij}^{\a}(x)\frac{\p^2 u_{\e}}{\p x_i \p x_j}-\sum_{i=1}^2 b_i^{\a}(x)\frac{\p u_{\e}}{\p x_i}\}=0
\end{equation}

$$
\qquad \hbox{in} \quad \Omega_{\e}=\{(x_1,x_2)|\quad -a\le x_1\le a,\quad 
f_0(x_1)+\e f_1(x_1,\frac{x_1}{\e})\le x_2\le b \}\subset \R^2,
$$
\begin{equation}\label{Gb}
<\n u_{\e},{\bf n}_{\e}>+c(x_1,\frac{x_1}{\e})u_{\e}=g(x_1,\frac{x_1}{\e}) 
\end{equation}
$$
\qquad \hbox{on}\quad \G_{\e}=\{(x_1,x_2)|\quad -a\le x_1\le a,\quad x_2=f_0(x_1)+\e f_1(x_1,\frac{x_1}{\e})\},
$$
\begin{equation}\label{Gd}
u_{\e}=0 \qquad \hbox{on}\quad \p\Omega_{\e}\backslash \G_{\e},
\end{equation}
where $\e>0$, $a_{ij}^{\a}(x)$, $b_i^{\a}(x)$ are Lipschitz in $x$ satisfying (\ref{lips}),  
  ${\bf n}_{\e}(x)$ is the outward unit normal to $\Omega_{\e}$, 
\begin{equation}\label{per}
c,\quad g,\quad f_1(x_1,\xi_1) \quad \hbox{are defined in} \quad \Omega_{\e}\times \R,
\quad \hbox{periodic in}\quad \xi_1\in \R\backslash {\bf Z},
\end{equation}
\begin{equation}\label{lc}
0\le f_1(x_1,\xi_1),\quad 0<C<c(x,\xi_1)  \qquad \hbox{in}\quad \Omega_{\e}\times \R\backslash {\bf Z}, 
\end{equation}
where $C>0$ is a constant, 
\begin{equation}\label{fd}
f_0'(\pm a)=0,\quad \frac{\p f_1}{\p \xi_1}(\pm a,\xi_1)=0,
\end{equation}
denoting $A_{\a}=(a_{ij}^{\a}(x))_{1\le i,j\le n}$, 
\begin{equation}\label{auni}
\l_1\le A_{\a}\le \L_1 \qquad \hbox{any}\quad  \a\in \CA.
\end{equation}

	We are interested in the limit of $u_{\e}$ of (\ref{Go})-(\ref{Gd}) as $\e$ goes 
to $0$. Remark that 
this problem is a straightforward generalization of Example 1.2, a similar case of which was treated in \cite{af} by the variational method. For our nonlinear problem, we need further assumptions 
listed in the following. These assumptions come from the formal asymptotic expansion 
 of $u_{\e}$ which we describe in below. (See also Remark 5.1 and Lemma 5.1 in below.)
 
\begin{equation}\label{ass1}
b_1^{\a}\equiv 0, \quad b_2^{\a}=a_{11}^{\a}f_0'' \qquad \hbox{any}\quad \a\in A, \quad 
x\in \Omega_{\e}, 
\end{equation} 
\begin{equation}\label{ass2}
\{a_{11}^{\a}(1+f_0'^2)-2a_{12}^{\a}f_0'+a_{22}^{\a}\}^2\ge 4(a_{11}^{\a}a_{22}^{\a}-{a_{12}^{\a}}^2)\qquad \hbox{forall}\quad \a\in A,\quad 
	x\in \Omega_{\e},
\end{equation}
 and for 
$$
O(x_1)=\{(\xi_1,\xi_2)|\quad \xi_2\ge f_1(x_1,\xi_1), \quad \hbox{periodic in} \quad \xi_1 \},
$$
\begin{equation}\label{ass4}
\p O(x_1) \quad {is} \quad C^{3,1}.
\end{equation}
	The existence and uniqueness of $u_{\e}$ ($\e>0$) is established in the general 
viscosity solutions theory. (See \cite{users}n.) Our goal is to show the existence of 
$u(x)$ such that 
\begin{equation}\label{conv}
\lim_{\e\to 0}u_{\e}(x)=u(x) \qquad \hbox{uniformly in}\quad \overline{\Omega},
\end{equation}
where $\Omega$$=\{(x_1,x_2)| \quad -a\le x_1\le a,\quad f_0 (x_1)\le x_2 \le b\}$, and to find the 
effective limit P.D.E. and B.C. for $u$. As for (\ref{conv}), we remark 
that our convergence is in $L^{\infty}$, while in \cite{af}n the convergence was in $H^1$. 
 The limit (effective) P.D.E. and B.C. are given  by using the 
 long time averaged result in $\S$ 3. Let us begin by deriving the cell problem for 
(\ref{Go})-(\ref{Gd}) by the 
 formal asymptotic expansions method:\\
\begin{equation}\label{formal}
u_{\e}=u(x)+\e v(\frac{x_1}{\e},\frac{x_2-f_0(x_1)}{\e})+O(\e^2),
\end{equation}
where we are assuming that "the corrector" $v$ depends only on $\xi_1=\frac{x_1}{\e}$ and $\xi_2=\frac{x_2-f_0(x_1)}{\e}$ ($\xi_1$, $\xi_2$ are rescaled variables.) From (\ref{formal}), we obtain 
\begin{eqnarray}\label{der1}
\frac{\p u_{\e}}{\p x_1}&=&\frac{\p u}{\p x_1}+\frac{\p v}{\p \xi_1}-f_0'(x_1)\frac{\p v}{\p \xi_2}+O(\e),\nonumber \\
\frac{\p u_{\e}}{\p x_2}&=&\frac{\p u}{\p x_2}+\frac{\p v}{\p \xi_2}+O(\e),
\end{eqnarray}
\begin{eqnarray}\label{der2}
\frac{\p^2 u_{\e}}{\p x_1^2}&=&\frac{\p^2 u}{\p x_1^2}-f_0''(x_1)\frac{\p v}{\p \xi_2}+\frac{1}{\e}\{\frac{\p^2 v}{\p \xi_1^2}-2f_0'(x_1)\frac{\p^2 v}{\p \xi_1 \p \xi_2}+(f_0')^2\frac{\p^2 v}{\p \xi_2^2}\}+
O(\e),\nonumber \\
\frac{\p^2 u_{\e}}{\p x_1 \p x_2}&=&\frac{\p^2 u}{\p x_1\p x_2}+\frac{1}{\e}(\frac{\p^2  v}{\p \xi_1\p \xi_2}-f_0'(x_1)\frac{\p^2 v}{\p \xi_2^2})+O(\e),\nonumber\\
\frac{\p^2 u_{\e}}{\p x_2^2}&=&\frac{\p^2 u}{\p x_2^2}+\frac{1}{\e} \frac{\p^2  v}{\p \xi_2^2} +O(\e),
\end{eqnarray}	
	First, by introducing (\ref{der1}) and (\ref{der2}) into 
\begin{eqnarray*}
&&-\sum_{i,j=1}^{2} a_{ij}^{\a} \frac{\p^2 u_{\e}}{\p x_i \p x_j}-\sum_{i=1}^2 b_i^{\a}
\frac{\p u_{\e}}{\p x_i}=\nonumber\\
&&=-\{
a_{11}^{\a}\frac{\p^2 u}{\p x_1^2}+2a_{12}^{\a}\frac{\p^2 u}{\p x_1 \p x_2}+a_{22}^{\a}\frac{\p^2 u}{\p x_2^2}-a_{11}^{\a}f_0''(x_1)\frac{\p v}{\p \xi_2}
\nonumber
\\
&& +b_1^{\a}(\frac{\p u}{\p x_1}+\frac{\p v}{\p \xi_1}-f_0'\frac{\p v}{\p \xi_2})+b_2(\frac{\p u}{\p x_2}+\frac{\p v}{\p \xi_2})
\} \nonumber\\
&&-\frac{1}{\e}[a_{11}^{\a}\{\frac{\p^2 v}{\p \xi_1^2}-2f_0'(x_1)\frac{\p^2 v}{\p \xi_1 \p \xi_2}+(f_0')^2 \frac{\p^2 v}{\p \xi_2^2}\}+2a_{12}^{\a}(\frac{\p^2 v}{\p \xi_1\p \xi_2}-f_0'(x_1)\frac{\p^2 v}{\p \xi_2^2})\nonumber\\
&&+a_{22}^{\a}\frac{\p^2 v}{\p \xi_2^2}
]\nonumber\\
\end{eqnarray*}
 and by using (\ref{ass1}),
\begin{eqnarray}\label{asym}
&&=-(a_{11}^{\a}\frac{\p^2 u}{\p x_1^2}+2a_{12}^{\a}\frac{\p^2 u}{\p x_1 \p x_2}+a_{22}^{\a}\frac{\p^2 u}{\p x_2^2})\\
&&-\frac{1}{\e}[a_{11}^{\a}
\frac{\p^2 v}{\p \xi_1^2}+2(a_{12}^{\a}-a_{11}^{\a}f_0'(x_1))
\frac{\p^2 v}{\p \xi_1 \p\xi_2}+\{a_{11}^{\a}(f_0')^2-2a_{12}^{\a}f_0'+a_{22}^{\a}\}\frac{\p^2 v}{\p\xi_2^2}]. \nonumber
\end{eqnarray}

{\bf Remark 5.1.} The condition (\ref{ass1}) was used to efface  the 
  dependence on $\xi$ (microscopic variable) in  the ordinary order ($O(1)$) part
in (\ref{asym}).\\

\bigskip  

	Let $(x,r,p)\in \Omega\times \R\times \R^2$ ($p=(p_1,p_2)$) be arbitrarily fixed, and define the following operators.
\begin{eqnarray}\label{opa}
&&P_{x,r,p}^{\a}(D^2_{\xi} v(\xi_1,\xi_2))\equiv \\
&&\equiv -[a_{11}^{\a}\frac{\p^2 v}{\p\xi_1^2} 
+2(a_{12}^{\a}-a_{11}^{\a}f_0')\frac{\p^2 v}{\p \xi_1 \p \xi_2} +
\{a_{11}^{\a}(f_0')^2-2a_{12}^{\a}f_0'+a_{22}^{\a}\}\frac{\p^2 v}{\p \xi_2^2}] \nonumber
\end{eqnarray}
in $O(x_1)$, and
\begin{equation}\label{ope}
P_{x,r,p}(D^2_{\xi} v(\xi_1,\xi_2))\equiv \sup_{\a\in \CA}\{P_{x,r,p}^{\a}(D^2_{\xi} v(\xi_1,\xi_2)\} \qquad \hbox{in} \quad O(x_1).
\end{equation}
	Next, by introducing (\ref{der1}) into (\ref{Gb}), we have 
\begin{eqnarray*}
&&\frac{1}{\sqrt{1+(f_0'+\frac{\p f}{\p \xi_1})^2}}
\{(f_0'+\frac{\p f_1}{\p \xi_1})\frac{\p u}{\p x_1}-\frac{\p u}{\p x_2}\}\nonumber\\
&&=g(x,\xi_1)-c(x,\xi_1)u-\frac{1}{\sqrt{1+(f_0'+\frac{\p f_1}{\p \xi_1})^2}}
\{(f_0'+\frac{\p f_1}{\p \xi_1})(\frac{\p v}{\p \xi_1}-f_0'\frac{\p v}{\p \xi_2})
-\frac{\p v}{\p \xi_2}\}. 
\end{eqnarray*}
By denoting the outward unit normal to the boundary of 
$$
\Omega=\{(x_1,x_2)|\quad -a\le x_1\le a,\quad x_2\ge f_0(x_1)\}
$$
 as 
$$
	\nu=\frac{1}{\sqrt{1+(f_0')^2}}(f_0',-1),
$$
the above equation on the boundary becomes 
\begin{eqnarray}\label{nub}
<\n u,\nu>&=&\frac{1}{\sqrt{1+(f_0')^2}}[-\frac{\p u}{\p x_1}\frac{\p f_1}{\p \xi_1}
-\sqrt{1+(f_0'+\frac{\p f_1}{\p \xi_1})^2}(cu-g)\nonumber\\
&&-(f_0'+\frac{\p f_1}{\p \xi_1})\frac{\p v}{\p \xi_1}+\{f_0'(f_0'+\frac{\p f_1}{\p \xi_1})+1\}\frac{\p v}{\p \xi_2}].
\end{eqnarray}
Let 
\begin{eqnarray}\label{defg}
\g(\xi_1,\xi_2)=\frac{(f_0'+\frac{\p f_1}{\p \xi_1},-\{f_0'(f_0'+\frac{\p f_1}{\p \xi_1})+1\})}{\sqrt{1+(f_0')^2}}\qquad \hbox{on}\quad \p O(x_1),
\end{eqnarray}
and for $(x,r,p)\in \Omega\times {\bf R}\times {\bf R}^2$
\begin{equation}\label{H}
H(x,r,p,\xi)=\frac{1}{\sqrt{1+(f_0')^2}}\{-\sqrt{1+(f_0'+\frac{\p f_1}{\p \xi_1})^2}(c(x,\xi_1)r-g)-p_1\frac{\p f_1}{\p \xi_1}\}.
\end{equation}
Then, (\ref{nub}) becomes 
\begin{eqnarray}\label{simple}
<\n u,\nu>=-\{<\g,\n_{\xi}v>-H(x,r,p,\xi)\}. 
\end{eqnarray}
From (\ref{asym}), (\ref{opa}), (\ref{ope}) and (\ref{simple}), the cell 
 problem for (\ref{Go})-(\ref{Gd}) should be the following: for any fixed $(x,r,p)\in \Omega \times {\bf R}\times {\bf R}^n$, find a unique number $d(x,p,r)$ such that 
 the following problem has a viscosity solution (corrector) $v(\xi_1,\xi_2)$.
\begin{eqnarray}\label{cell}
&&P_{x,r,p}(D^2_{\xi}v(\xi_1,\xi_2))=0 \qquad \hbox{in} \quad O(x_1),\nonumber\\
&&d(x,r,p)+<\n_{\xi} v,\g>-H(x,r,p,\xi)=0 \qquad \hbox{on}\quad \p O(x_1),\nonumber\\
&&v \quad \hbox{is bounded in}\quad \overline{O(x_1)}.
\end{eqnarray}

\begin{lemma}{\bf Lemma 5.1.} Let (\ref{ass2}) hold. Then, the operators $P_{x,r,p}^{\a}(\xi_1,\xi_2)$ are uniformly elliptic operators uniformly in $\a\in A$: 
there exist constants $0< \l_1'<\L_1'$ such that 
$$
\l_1' I\le
\left(
\begin{array}{cc}
a_{11}^{\a} & a_{12}^{\a}-a_{11}^{\a}f_0'\\
a_{12}^{\a}-a_{11}^{\a}f_0' & a_{22}^{\a}-2a_{12}^{\a}f_0'+a_{11}^{\a}f_0'
\end{array}
\right)
\le \L_1' I \qquad \hbox{any}\quad \a\in A.
$$
\end{lemma}

{\it Proof of Lemma 5.1.} The claim can easily confirmed by an elementary calculation. 
And we leave it to the readers.\\

\begin{lemma}{\bf Lemma 5.2.}
 Let $\a\in \CA$ and $(x,r,p)$ be fixed, and let $O(x_1)$, $P_{x,r,p}^{\a}(D^2_{\xi})$, $\g(\xi)$ and $H(x,r,p,\xi)$ be defined in (\ref{ass4}), (\ref{opa}), (\ref{defg}) and (\ref{H}).  
Assume that (\ref{per})-(\ref{ass4}) hold. Then, there exists a unique number $d^{\a}(x,r,p)$ 
such that the following problem has a viscosity solution $v(\xi_1,\xi_2)$.
\begin{eqnarray}\label{cella}
&&P_{x,r,p}^{\a}(D^2_{\xi}v(\xi_1,\xi_2)=0 \qquad \hbox{in} \quad O(x_1),\nonumber\\
&&d^{\a}(x,r,p)+<\n_{\xi} v,\g>-H(x,r,p,\xi)=0 \qquad \hbox{on}\quad \p O(x_1),\nonumber\\
&&v\quad\hbox{is bounded in}\quad \overline{O(x_1)}.
\end{eqnarray}
\end{lemma}

{\it Proof of Lemma 5.2.} From (\ref{defg}), we confirm easily that there exists a positive constant 
 $\g_1>0$ such that 
$$
<\g,\zeta>\quad >\g_1>0 \qquad \hbox{on}\quad \p O(x_1),
$$
where $\zeta=\frac{(\frac{\p f_1}{\p \xi_1},-1)}{\sqrt{(\frac{\p f_1}{\p \xi_1})^2+1}}$ 
the outward unit normal to $\p O(x_1)$. Then from Theorem 3.3, there exists a unique number 
$d^{\a}(x,r,p)$ such that (\ref{cella}) has a  viscosity solution $v$. 

\begin{lemma}{\bf Lemma 5.3.} We assume the same assumptions as in Lemma 5.2. For any fixed $(x,r,p)$, there exists a unique number $d$ such that (\ref{cell}) has a viscosity solution $v(\xi_1,\xi_2)$.
 Moreover,
\begin{equation}\label{as}
d(x,r,p)\le d^{\a}(x,r,p) \qquad \hbox{any} \quad \a\in \CA.
\end{equation}
\end{lemma}
{\it Proof of Lemma 5.3.} From Theorem 3.3, there exists a unique number $d(x,r,p)$ such that 
 (\ref{cell}) has a viscosity solution $v$. 
 The inequality (\ref{as}) comes from the construction of the number $d$ and $d^{\a}$ in 
the proofs of Proposition 3.2. and Theorem 3.3. That is, 
$$
	d=\lim_{R\to \infty} d_R, \quad d^{\a}=\lim_{R\to \infty} d^{\a}_R,
$$ 
where $d$ and $d_R$ ($R\in {\bf N}$) are characterized by the following: for 
$O_R(x_1)=O(x_1)\cap \{\xi_2\le R\}$ 
\begin{eqnarray}
&&P_{x,r,p} (D^2_{\xi}v_R(\xi_1,\xi_2)=0 \qquad \hbox{in} \quad O_R(x_1),\nonumber\\
&&d_R(x,r,p)+<\n_{\xi} v_R,\g>-H(x,r,p,\xi)=0 \qquad \hbox{on}\quad \p O(x_1),\nonumber\\
&&<\n_{\xi} v_R,{\bf n}>=0\qquad \hbox{on}\quad \{\xi_2=R\},\nonumber
\end{eqnarray}
and
\begin{eqnarray}
&&P_{x,r,p}^{\a}(D^2_{\xi}v_R^{\a}(\xi_1,\xi_2)=0 \qquad \hbox{in} \quad O_R(x_1),\nonumber\\
&&d^{\a}_R(x,r,p)+<\n_{\xi} v_R^{\a},\g>-H(x,r,p,\xi)=0 \qquad \hbox{on}\quad \p O(x_1),\nonumber\\
&&<\n_{\xi} v_R^{\a},{\bf n}>=0\qquad \hbox{on}\quad \{\xi_2=R\},\nonumber
\end{eqnarray}
where $\bf n$ is the outward unit normal to $\p O_R(x_1)$ on $\{\xi_2=R\}$. From the stochastic 
representations (\ref{long}) of $d_R$ and $d^{\a}_R$ in the approximating problems  (\ref{pR}), we see that 
$$
d_R \le d^{\a}_R \qquad \hbox{any}\quad R\in {\bf N}. 
$$
 Therefore, (\ref{as}) was proved.\\

	Since the oscillating Neumann boundary condition prevent us from obtaining  the uniform 
gradient bounds of $u_{\e}$ ($\e>0$), we need to treat the upper and lower envelopes. 

\begin{lemma}{\bf Lemma 5.4.} Assume that (\ref{lips}), (\ref{per})-(\ref{ass4}) hold. Let $u_{\e}$ be the solution of (\ref{Go})-(\ref{Gd}). Then, there exists a constant $M>0$ such that 
\begin{equation}\label{Mo}
|u_{\e}|<M \qquad \hbox{any} \quad \e\in (0,1).
\end{equation}
\end{lemma}

{\it Proof of Lemma 5.4.} Let $x_0=(0,b+r)\in \R^2$, where $r>0$. Define
$$
v(x)=A(r^{-p}-|x-x_0|^{-p}) \qquad x\in \Omega_{\e}.
$$
Then, for $A>0$ large enough, $v$ is a super solution of (\ref{Go})-(\ref{Gd}) for any $\e\in (0,1)$. From the comparison result for (\ref{Go})-(\ref{Gd}), we get (\ref{Mo}).\\

From (\ref{Mo}), 
$$
u^{\ast}(x)=\limsup_{\e\downarrow 0,y\to x} u_{\e}(y), \quad
u_{\ast}(x)=\liminf_{\e\downarrow 0,y\to x} u_{\e}(y) \qquad x\in \overline{\Omega},
$$
are well-definded. Moreover, from (\ref{auni}) and the Krylov-Safonov inequality we can extract 
a subsequence ${\e}'\to 0$ such that 
\begin{equation}\label{match}
\lim_{{\e}'\downarrow 0} u_{{\e}'}=u \qquad \hbox{locally uniformly in} \quad \Omega,\quad 
u^{\ast}\ge u \ge u_{\ast}.
\end{equation}

We claim the following.
\begin{lemma}{\bf Lemma 5.5.} Assume that (\ref{per})-(\ref{ass4}) hold. Then, $u^{\ast}$ and 
$u_{\ast}$ are respectively viscosity sub and super solutions of the following problem. 
\begin{equation}\label{elim}
\sup_{\a\in \CA} \{-\sum_{i,j=1}^n a_{ij}^{\a}\frac{\p^2 u}{\p x_i \p x_j}-
\sum_{i=1}^n b_i^{\a} \frac{\p u}{\p x_i}\}=0 \qquad \hbox{in}\quad \Omega,
\end{equation}
\begin{equation}\label{blim}
<\n u, \nu>+\overline{L}(x,u,\n u)=0 \qquad \hbox{on} \quad \G_0,
\end{equation}
where $\nu$ is the outward unit normal to $\Omega$ defined on
$$
\G_0=\{(x_1,x_2)|\quad -a\le x_1\le a,\quad x_2=f_0(x_1)\},
$$
and for $(x,r,p)\in \overline{\Omega}\times \R \times \R^2$,
\begin{equation}\label{Lo}
\overline{L}(x,r,p)=-d(x,r,p),
\end{equation}
where $d(x,r,p)$ is defined in (\ref{cell}).
\end{lemma}

{\it Proof of Lemma 5.5.} From (\ref{match}) and by the usual stability results of the 
viscosity solutions, it is clear that (\ref{elim}) holds. In the following, we shall  see (\ref{blim}).\\

\underline{Step 1.} We shall show that $u^{\ast}$ satisfies 
$$
<\n u^{\ast},\nu>+\overline{L}(x,\n u^{\ast},\H u^{\ast})\le 0 \qquad \hbox{on}\quad 
\G_0,
$$
in the sense of viscosity solutions. Remark that $\Omega_{\e}\subset \Omega$ for any 
$\e\in [0,1)$. 
 Let $\phi\in C^2(\overline{\Omega})$ be such that $u^{\ast}-\phi$ takes its strict maximum at 
$x_0=(x_{01},x_{02})$$\in \G_0$ with  $u^{\ast}(x_0)=\phi(x_0)$. From the  definition 
of the  Neumann type boundary condition in the sense of viscosity solutions, we are to show either 
\begin{equation}\label{p1}
\sup_{\a\in \CA}\{-\sum_{ij}a_{ij}^{\a}\frac{\p^2 \phi}{\p x_i \p x_j}(x_0)-
\sum_{i}b_i^{\a}\frac{\p \phi}{\p x_i}(x_0)\}\le 0,
\end{equation} 
or
\begin{equation}\label{p2}
<\n \phi(x_0),\nu>+\overline{L}(x_0,\n \phi(x_0),\H \phi(x_0))\le 0.
\end{equation}
	We shall assume that both (\ref{p1}) and (\ref{p2}) are not true, and  shall  seek a 
contradiction. Thus, assume there exist constants $\theta_1$ and $\theta_2$ such that 
\begin{equation}\label{ct1}
\sup_{\a\in \CA}\{-\sum_{ij}a_{ij}^{\a}\frac{\p^2 \phi}{\p x_i \p x_j}(x_0)-
\sum_{i}b_i^{\a}\frac{\p \phi}{\p x_i}(x_0)\}\equiv \theta_1>0,
\end{equation}
\begin{equation}\label{ct2}
<\n \phi(x_0),\nu>+\overline{L}(x_0,\n \phi(x_0),\H \phi(x_0))\equiv \theta_2>0.
\end{equation}

 For $(x_0,r_0,p_0)=(x_0,\phi(x_0),\n \phi(x_0))$, from Lemma 5.2 there exists a number $d(x_0,r_0,p_0)$ and $v$ of
\begin{eqnarray}\label{A}
&&P_{x_0,r_0,p_0}(D^2_{\xi} v(\xi_1,\xi_2))=0 \qquad \hbox{in} \quad O(x_{01}),\\
&&d(x_0,r_0,p_0)+<\n_{\xi} v,\g>-H(x_0,r_0,p_0,\xi)=0 \qquad \hbox{on}\quad \p O(x_1).\nonumber
\end{eqnarray}
Since $\xi_2\ge f_1 (x_1,\xi_1)$ for any $(\xi_1,\xi_2)\in O(x_1)$, we may define 
$$
\phi_{\e}(x_1,x_2)=\phi(x_1,x_2)+\e v(\frac{x_1}{\e},\frac{x_2-f_0(x_1)}{\e}) \qquad \hbox{in}\quad 
 \overline{\Omega_{\e}}.
$$
 We claim that $\phi_{\e}$ is the viscosity supersolution of 
\begin{equation}\label{but1}
\sup_{\a\in \CA}\{-\sum_{ij}a_{ij}^{\a}\frac{\p^2 \phi_{\e}}{\p x_i \p x_j}-
\sum_{i}b_i^{\a}\frac{\p \phi_{\e}}{\p x_i}\}>\frac{1}{4}\theta_1 \qquad 
\hbox{in}\quad B(x_0,r)\cap \Omega_{\e},
\end{equation}
\begin{equation}\label{but2}
<\n \phi_{\e},{\bf n}_{\e}>+c(x,\frac{x_1}{\e})\phi_{\e}
-g(x,\frac{x_1}{\e})>\frac{1}{4}\theta_2
\qquad \hbox{on}\quad B(x_0,r)\cap \G_{\e},
\end{equation}
in the sense of viscosity solutions in some small neighborhood of $x_0$, $B(x_0,r)$ ($r>0$ is uniform in $\e\in (0,1)$). To see this, assume for $\psi\in C^2(\overline{\Omega})$, $\phi_{\e}-\psi$ takes its 
minimum at $(\overline{x_1},\overline{x_2})$ with 
$\phi_{\e}(\overline{x_1},\overline{x_2})=\psi_{\e}(\overline{x_1},\overline{x_2})$. \\

First, let us assume that $(\overline{x_1},\overline{x_2})\in \Omega_{\e}$. 
 We write 
\begin{equation}\label{ye}
\eta(\xi_1,\xi_2)\equiv \frac{1}{\e}(\psi-\phi)(\e \xi_1,\e \xi_2+f_0(\e \xi_1)) 
\qquad (\xi_1,\xi_2)\in O(x_1),
\end{equation}
$$
\overline{\xi_1}\equiv 
\frac{\overline{x_1}}{\e},\quad  \overline{\xi_2}\equiv \frac{\overline{x_2}-f_0(\overline{x_1})}{\e}. 
$$
Hence, 
$$
(v-\eta)(\overline{\xi_1},\overline{\xi_2})\le (v-\eta)(\xi_1,\xi_2),
$$
in a neighborhood of $(\frac{x_{01}}{\e},\frac{x_{02}-f_0(x_{01})}{\e})$
$\equiv (\xi_{01},\xi_{02})$. 
 Now, from (\ref{ye}),
\begin{eqnarray}\label{dif1}
\frac{\p \eta}{\p \xi_1}&=&\frac{\p}{\p x_1}(\psi-\phi)+\frac{\p}{\p x_2}(\psi-\phi)f_0'(\e \xi_1),\nonumber \\ 
\frac{\p \eta}{\p \xi_2}&=&\frac{\p}{\p x_2}(\psi-\phi),
\end{eqnarray}
\begin{eqnarray}\label{dif2}
\frac{\p^2 \eta}{\p \xi_1^2}&=&\e\{
\frac{\p^2}{\p x_1^2}(\psi-\phi)+2\frac{\p^2}{\p x_1 \p x_2}(\psi-\phi)f_0'+
\frac{\p^2}{\p x_2^2}(\psi-\phi)(f_0')^2 \nonumber\\
&+&\frac{\p}{\p x_2}(\psi-\phi)f_0''
\}, \nonumber\\
\frac{\p^2 \eta}{\p \xi_1 \p \xi_2}&=&\e\{
\frac{\p^2}{\p x_1 \p x_2}(\psi-\phi)+
\frac{\p^2}{\p x_2^2}(\psi-\phi)(f_0')
\},\nonumber\\
\frac{\p^2 \eta}{\p \xi_2^2}&=&\e
\frac{\p^2}{\p x_2^2}(\psi-\phi). 
\end{eqnarray}
 Since $v(\xi_1, \xi_2)$ is the viscosity solution of (\ref{A}), by (\ref{ye}),  (\ref{dif1}) and (\ref{dif2}), for any $\d>0$ there exists a control $\overline{\a}\in A$ 
 such that 
\begin{eqnarray*}
-[a_{11}^{\overline{\a}}\{ \frac{\p^2}{\p x_1^2}(\psi-\phi)&+&2\frac{\p^2}{\p x_1 \p x_2}(\psi-\phi)f_0'+
\frac{\p^2}{\p x_2^2}(\psi-\phi)(f_0')^2+\frac{\p}{\p x_2}(\psi-\phi)f_0''\}\\
&&+2(a_{12}^{\overline{\a}}-a_{11}^{\overline{\a}}f_0')\{ \frac{\p^2}{\p x_1 \p x_2}(\psi-\phi)+
\frac{\p^2}{\p x_2^2}(\psi-\phi)(f_0')\}\\
&&+ (a_{22}^{\overline{\a}}-2a_{12}^{\overline{\a}}f_0'+a_{11}^{\overline{\a}}(f_0')^2)
\frac{\p^2}{\p x_2^2}(\psi-\phi)(\overline{x_1},\overline{x_2})]\ge -\d.
\end{eqnarray*}
We can simplify the above by using $a_{11}^{\overline{\a}}f_0''=b_2^{\overline{\a}}$ ((\ref{ass1})) to
\begin{eqnarray}
&&(-\sum_{ij}a_{ij}^{\overline{\a}}(x_0) \frac{\p^2 \psi}{\p x_i \p x_j}-
\sum_{i}b_{i}^{\overline{\a}}(x_0)\frac{\p \psi}{\p x_i}+
\sum_{ij}a_{ij}^{\overline{\a}}(x_0) \frac{\p^2 \phi}{\p x_i \p x_j}\nonumber\\
&&+
\sum_{i}b_{i}^{\overline{\a}}(x_0)\frac{\p \phi}{\p x_i})(\overline{x_1},\overline{x_2})\ge -\d.\nonumber
\end{eqnarray}
Thus, since $\d>0$ is arbitrary, 
\begin{eqnarray}
 &&\sup_{\a\in \CA}\{
-\sum_{ij}a_{ij}^{\a}(x_0) \frac{\p^2 \psi}{\p x_i \p x_j}-
\sum_{i}b_{i}^{\a}(x_0)\frac{\p \psi}{\p x_i}\}(\overline{x_1},\overline{x_2})
\nonumber\\
&&\ge 
(-\sum_{ij}a_{ij}^{\overline{\a}}(x_0) \frac{\p^2 \psi}{\p x_i \p x_j}-
\sum_{i}b_{i}^{\overline{\a}}(x_0)\frac{\p \psi}{\p x_i})(\overline{x_1},\overline{x_2}) 
\nonumber\\
&&\ge -\d +
(-\sum_{ij}a_{ij}^{\overline{\a}}(x_0) \frac{\p^2 \phi}{\p x_i \p x_j}-
\sum_{i}b_{i}^{\overline{\a}}(x_0)\frac{\p \phi}{\p x_i})(\overline{x_1},\overline{x_2})\ge \frac{\theta_1}{2},\nonumber
\end{eqnarray}
for $(\overline{x_1},\overline{x_2})$ is near to $x_0$, and for  $r>0$ small enough. 
Therefore, (\ref{but1}) was shown. \\

	Next, we assume 
\begin{equation}\label{ong}
(\overline{x_1},\overline{x_2})\in \G_{\e}.
\end{equation}
 Again, we use the same function $\eta$ defined in (\ref{ye}) and denote  
$\xi_1=\frac{x_1}{\e}$, $\xi_2=\frac{x_2-f_0(x_1)}{\e}$, 
$$
(\overline{\xi_1},\overline{\xi_2})=
(\frac{\overline{x_1}}{\e},\frac{\overline{x_2}-f_0(\overline{x_1})}{\e}),  \quad (\overline{\xi_{01}}, \overline{\xi_{02}})=
(\frac{\overline{x_{01}}}{\e},\frac{\overline{x_{02}}-f_0(\overline{x_{01}})}{\e}). 
$$ 
 Thus,
\begin{equation}\label{nbd}
(v-\eta)(\overline{\xi_1},\overline{\xi_2})\le (v-\eta)(\xi_1,\xi_2),
\end{equation}
in a small neighborhood of $(\overline{\xi_{01}},\overline{\xi_{02}})$. By (\ref{ong}) 
$\overline{x_2}=f_0(\overline{x_1})+\e f_1(\overline{x},\frac{\overline{x_1}}{\e})$, and 
$$
\overline{\xi_2}=f_1(\overline{x},\overline{\xi_1}), \quad (\overline{\xi_1},\overline{\xi_2})\in \p O(x_{1}).
$$
 Since $v$ satisfies (\ref{A}), from the 
 definition of the viscosity solution 
\begin{equation}\label{super1}
P_{x_0,\phi(x_0),\n \phi(x_0)} (D^2_{\xi} \eta)(\overline{\xi_1},\overline{\xi_2})\ge 0,
\end{equation} 
or
\begin{equation}\label{super2}
d(x_0,\phi (x_0),\n \phi (x_0))+<\n_{\xi} \eta,\g>(\overline{\xi_1},\overline{\xi_2})
-H(x_0,\phi(x_0),\n \phi(x_0),\overline{\xi_1},\overline{\xi_2})\ge 0.
\end{equation}
 In the case of (\ref{super1}), as before we obtain 
\begin{equation}\label{same}
\sup_{\a\in A} \{
-\sum_{ij} a_{ij}^{\a}(\overline{x})\frac{\p^2 \psi}{\p x_i \p x_j}(\overline{x})
-\sum_{i} b_{i}^{\a}(\overline{x})\frac{\p \psi}{\p x_i}(\overline{x})
\}>\frac{1}{4}{\theta}_1.
\end{equation}
 In the case of (\ref{super2}),   from (\ref{Lo}), (\ref{H}) and (\ref{super2}),
\begin{equation}\label{paris}
-\overline{L}(x_0,\phi(x_0),\n \phi(x_0))+\frac{1}{\sqrt{(f_0')^2+1}}
<\n_{\xi}\eta,(f_0'+\frac{\p f_1}{\p \xi_1},
-f_0'(f_0'+\frac{\p f_1}{\p \xi_1})-1)>
\end{equation}
$$
-\frac{1}{\sqrt{(f_0')^2+1}}(-\sqrt{1+(f_0'+\frac{\p f_1}{\p \xi_1})^2}c(x,\xi_1)\phi 
-\frac{\p\phi}{\p x_1}\frac{\p f_1}{\p \xi_1}+ \sqrt{1+(f_0'+\frac{\p f_1}{\p \xi_1})^2}g)
\ge 0.
$$
 Introducing (\ref{dif1}) to (\ref{paris}) 
$$
-\overline{L}(x_0,\phi(x_0),\n \phi(x_0))+\frac{1}{\sqrt{(f_0')^2+1}}
<\n (\psi -\phi)(x_0),(f_0'+\frac{\p f_1}{\p \xi_1},-1)>
$$
$$
-\frac{1}{\sqrt{(f_0')^2+1}}(-\sqrt{1+(f_0'+\frac{\p f_1}{\p \xi_1})^2}c\phi-
\frac{\p \phi}{\p x_1}\frac{\p f_1}{\p \xi_1}+\sqrt{1+(f_0'+\frac{\p f_1}{\p \xi_1})^2 }g)\ge o(\e),
$$ and 
deviding the both hands sides of the above by $\sqrt{1+(f_0'+\frac{\p f_1}{\p \xi_1})^2}$, 
 by remarking that 
$$
{\bf n}_{\e}=(\frac{f_0'+\frac{\p f_1}{\p \xi_1}}{\sqrt{1+(f_0'+\frac{\p f_1}{\p \xi_1})^2}},
\frac{-1}{\sqrt{1+(f_0'+\frac{\p f_1}{\p \xi_1})^2}})+o(\e),
$$
we have
$$
\frac{1}{\sqrt{(f_0')^2+1}}
<\n \psi (x_0),{\bf n}_{\e}>-
\frac{\overline{L}(x_0,\phi(x_0),\n \phi(x_0))}{\sqrt{1+(f_0'+\frac{\p f_1}{\p \xi_1})^2 }}
\qquad \qquad \qquad \qquad \qquad \qquad \qquad \qquad $$
$$
\ge \frac{1}{\sqrt{1+(f_0'+\frac{\p f_1}{\p \xi_1})^2} }<\n \phi,\nu>
-\frac{1}{\sqrt{(f_0')^2+1}}c \phi + \frac{1}{\sqrt{(f_0')^2+1}} g+o(\e).
$$
 By using (\ref{ct2}) and  multiplying the both hands sides of the above by 
$\sqrt{(f_0')^2+1}$, we get 
$$
<\n \psi(x_0),{\bf n}_{\e}>+c\phi(x_0)-g\ge \overline{L}(x_0,\phi(x_0),\n \phi(x_0))+<\n \phi(x_0),\nu>
\equiv \theta_2>0,
$$
and for $r>0$ and $\e>0$ small enough,
\begin{equation}\label{toto}
<\n \psi(x_1),{\bf n}_{\e}>+c\phi(x_1)-g\ge \frac{1}{2}{\theta}_2.
\end{equation}
We have proved (\ref{but2}). Thus, in $B(x_0,r)\cap \overline{\Omega_{\e}}$, we have 
 (\ref{but1})-(\ref{but2}) and (\ref{1})-(\ref{2}). Therefore, 
$$
\max_{\overline{B(x_0,r)\cap \Omega_{\e}}} (u_{\e}-\phi_{\e})=
\max_{\p (B(x_0,r)\cap \Omega_{\e})} (u_{\e}-\phi_{\e}).
$$
From (\ref{Gb}) and 
(\ref{but2}), by using a similar argument in the proof of Lemma 2.6, 
$$
<\n (u_{\e}-\phi_{\e}),{\bf n}_{\e}>+c(u_{\e}-\phi_{\e})<-\frac{1}{4}\theta_2<0 \qquad 
\hbox{on} \quad \G_{\e}\cap B(x_0,r),
$$
in the sense of viscosity solutions. 
By letting $\e$ tends to zero, $\max(u_{\e}-\phi_{\e})$ goes to zero  and there exists $\e_0>0$ such that 
$$
<\n (u_{\e}-\phi_{\e}),{\bf n}_{\e}>\quad <-\frac{1}{8}\theta_2<0 \qquad \hbox{on}\quad \G_{\e}\cap B(x_0,r) \qquad \hbox{any}\quad \e\in (0,\e_0).
$$

 From this, if $u_{\e}-\phi_{\e}$ ($\e\in (0,\e_0)$) takes its local  maximum on $\G_{\e}\cap B_r(x_0)$ the strong maximum principle (Lemma A) leads a contradiction.
 Thus, $u_{\e}-\phi_{\e}$ must take its maximum on $\p\overline{B(x_0,r)\cap \Omega_{\e}}\backslash \Gamma_{\e}$, that is on $\p B(x_0,r)$. However this contradicts 
to the fact that $u-\phi$ takes its strong maximum in $\overline{B(x_0,r)\cap \Omega}$ 
at $x_0$. Thus, we proved (\ref{p1})-(\ref{p2}).\\

\underline{Step 2.} The fact that $u_{\ast}$ is a supersolution of
$$
<\n u_{\ast},\nu>+\overline{L}(x,\n u_{\ast},\H u_{\ast})\le 0 \qquad \hbox{on}\quad 
\G_0,
$$
in the sense of viscosity solutions can be shown similarly to (and slightly easier than)  Step 1. We omit the details, since the argument is parallel. \\

From the above, we complete the proof of Lemma 5.5.\\

\begin{lemma}{\bf Lemma 5.6.} Assume that (\ref{per})-(\ref{ass4}) hold. Then,
$$
u^{\ast}=u_{\ast}=0 \qquad x\in \p\Omega\backslash \G_0.
$$
\end{lemma}
{\it Proof of Lemma 5.6.} Let $x_0\in \p\Omega_{\e}\backslash \G_{\e}$ be arbitrarily fixed. We can take $\underline{v}$ and 
$\overline{v}$,  sub and super solutions of 
$$
\sup_{\a\in \CA}\{-\sum_{ij}a_{ij}^{\a}\frac{\p^2 \underline{v}}{\p x_i \p x_j}
-\sum_{i} b_i^{\a}\frac{\p \underline{v}}{\p x_i}\}\le 0 \qquad \hbox{in} \quad 
\Omega_{\e},
$$ 
$$
<\n \underline{v},{\bf n}_{\e}>+c\underline{v}\le g \qquad \hbox{on}\quad \G_{\e},
$$
$$
\underline{v}(x_0)=0,\quad \underline{v}(x)\le 0 \qquad \hbox{on}\quad 
\p\Omega\backslash \G_{\e}, 
$$
and
$$
\sup_{\a\in \CA}\{-\sum_{ij}a_{ij}^{\a}\frac{\p^2 \overline{v}}{\p x_i \p x_j}
-\sum_{i} b_i^{\a}\frac{\p \overline{v}}{\p x_i}\}\ge 0 \qquad \hbox{in} \quad 
\Omega_{\e},
$$ 
$$
<\n \overline{v},{\bf n}_{\e}>+c\overline{v}\ge g \qquad \hbox{on}\quad \G_{\e},
$$
$$
\overline{v}(x_0)=0,\quad \overline{v}(x)\ge 0 \qquad \hbox{on}\quad 
\p\Omega\backslash \G_{\e}. 
$$

 From the comparison, 
$$
\underline{v}\le u_{\e}\le \overline{v} \qquad \hbox{any} \quad \e\in (0,1),
$$
and thus
$$
\underline{v}\le u_{\ast}\le u^{\ast}\le \overline{v} \qquad \hbox{any}\quad  x\in \overline{\Omega}.
$$
In particular, at $x_0$,
$$
\underline{v}(x_0)=u_{\ast}(x_0)= u^{\ast}(x_0)= \overline{v}(x_0)=0.
$$

\begin{lemma}{\bf Lemma 5.7.}
 The function $\overline{L}(x,r,p)$ is increasing in $r$.
\end{lemma}

{\it Proof of Lemma 5.7.} From the definition of $\overline{L}$, we are to show that 
 $d(x,r,p)$ is decreasing in $r$. As we mentioned in the proof of (\ref{as}) in Lemma 5.3, 
this fact is clear from the construction of $d$ and its meaning in (\ref{long}).\\ 

From Lemmas 5.5-5.7, we arrive at the following result.

\begin{theorem}{\bf Theorem 5.8.} Assume that (\ref{per})-(\ref{ass4}) hold. Then, there exists a unique 
function $u(x)$ such that 
$$
\lim_{\e\downarrow 0} u_{\e}(x)=u(x) \qquad \hbox{locally uniformly in} \quad \overline{\Omega},
$$
which is the unique solution of (\ref{elim}), (\ref{blim}), and (\ref{Gd}).
\end{theorem}

{\it Proof of Theorem 5.8.} From Lemmas 5.5, 5.6 and 5.7, the limit $u^{\ast}=u_{\ast}=u$ is unique and is a  solution of 
the above problem. Moreover, since from Lemma 5.7 the uniqueness holds for (\ref{elim})-(\ref{blim}) and (\ref{Gd}), $u$  is the unique solution.  (We refer the readers to \cite{users}) and G. Barles \cite{barles} for such uniqueness results. And,  we proved the claim.\\

{\bf Remark 5.2.} The effective boundary condition (\ref{blim}) is in general nonlinear. 
 However, for the linear problem as in Example 1.2, (\ref{blim}) is lenear and matchs to the 
 result in \cite{af}.\\

\begin{example}{\bf Example 5.1.} Let $f_0'\equiv 0$, and assume that $a_{11}=a_{22}=1$, $a_{12}=0$. Then, 
$$
	\overline{L}(x,r,p)=-d(x,r,p),
$$ 
 is obtained by the following long time averaged problem: 
\begin{eqnarray*}
&&P_{x,r,p}(D^2_{\xi}v(\xi_1,\xi_2)) =-\frac{\p^2 v}{\p \xi_1^2}-
 \frac{\p^2 v}{\p \xi_2^2}=0 \qquad \hbox{in}\quad O(x_{1}),\\
&&d(x,r,p)-<\n_{\xi} v,(\frac{\p f_1}{\p \xi_1},-1)> -\{-\sqrt{1+(\frac{\p f_1}{\p \xi_1})^2}(c(x,\xi_1)r-g)-p_1 \frac{\p f_1}{\p \xi_1}\}=0
\end{eqnarray*}
$$
\qquad\qquad\qquad\qquad\qquad\quad\qquad\qquad\qquad\qquad\qquad\qquad\quad\qquad\qquad \hbox{in}\quad O(x_{1}),
$$
where 
$$
O(x_{1})=\{(\xi_1,\xi_2)|\quad \hbox{periodic in}\quad \xi_1\in {\bf R}\backslash{\bf Z},\quad \xi_2\ge f_1(x,\xi_1)\}.
$$
By integrating the above problem in $\xi_1\in [0,1]$, and by remarking that $f_1$ and $v$ are periodic in 
$\xi_1$, we have 
$$
d(x,r,p)=-r\int_0^1 \sqrt{1+(\frac{\p f_1}{\p \xi_1})^2}c(x,\xi_1) d\xi_1+\int_0^1 \sqrt{1+(\frac{\p f_1}{\p \xi_1})^2}g d\xi_1.
$$
 Therefore, 
$\overline{L}(x,r,p)$ is linear in $r$. 
\end{example}

{\bf Remark 5.3.} Although in this paper we considered a particular exaple of the oscillating Neumann condition  ((\ref{Gb})) in ${\bf R}^2$, we can apply the same method to more general  homogenization of the oscillating  boundary conditions in ${\bf R}^n$. We shall give more general formulation of this kind of problem in the future occassion.\\
 

\end{document}